\documentclass[toc=bib,toc=idx,abstract=false,parskip=half]{scrartcl}
\pdfoutput=1
\usepackage[utf8]{inputenc}
\usepackage[T1]{fontenc}
\usepackage[english]{babel}
\usepackage{csquotes}
\usepackage{mathrsfs}
\usepackage{amsmath}
\usepackage{amsthm}
\usepackage{thmtools}
\usepackage{amsfonts}
\usepackage{amssymb}
\usepackage{stmaryrd}
\usepackage{mathtools}
\usepackage{mathdots}
\usepackage{enumitem}
\usepackage{booktabs}
\usepackage{calc}
\usepackage{graphicx}
\usepackage[dvipsnames]{xcolor}
\usepackage[style=alphabetic,citestyle=alphabetic,maxcitenames=4,maxalphanames=4,maxbibnames=10,backref=false,hyperref=true,backend=biber,doi=false,isbn=false,giveninits=true]{biblatex}
\usepackage[colorlinks=true,linkcolor=black,urlcolor=black,citecolor=black,bookmarksnumbered]{hyperref}

\usepackage{tikz}
\usetikzlibrary{arrows,decorations.markings,chains,calc,matrix,patterns}

\bibliography{balancedideals}
\DeclareFieldFormat{postnote}{#1}	
\DeclareFieldFormat{multipostnote}{#1}
\DeclareFieldFormat[article,eprint,unpublished,incollection]{title}{\mkbibemph{#1}}
\DeclareFieldFormat[article,eprint,unpublished,incollection]{journaltitle}{#1}
\DeclareFieldFormat[incollection]{booktitle}{#1}

\makeindex

\AtBeginDocument{%

}

\declaretheoremstyle[spaceabove=\topsep,spacebelow=0pt,bodyfont=\normalfont]{scdef}
\declaretheoremstyle[spaceabove=\topsep,spacebelow=0pt,bodyfont=\itshape]{scthm}
\declaretheoremstyle[spaceabove=\topsep,spacebelow=0pt,headfont=\normalfont\itshape,notefont=\normalfont\itshape,notebraces={}{},headformat={\NAME\NOTE},postheadspace=1em,qed=\qedsymbol]{scprf}

\declaretheorem[style=scthm,numberwithin=section,name=Theorem,    refname={Theorem,Theorems},        Refname={Theorem,Theorems}]        {Thm}
\declaretheorem[style=scthm,sharenumber=Thm,     name=Lemma,      refname={Lemma,Lemmas},            Refname={Lemma,Lemmas}]            {Lem}
\declaretheorem[style=scthm,sharenumber=Thm,     name=Corollary,  refname={Corollary,Corollaries},   Refname={Corollary,Corollaries}]   {Cor}
\declaretheorem[style=scthm,sharenumber=Thm,     name=Proposition,refname={Proposition,Propositions},Refname={Proposition,Propositions}]{Prop}
\declaretheorem[style=scdef,sharenumber=Thm,     name=Definition, refname={Definition,Definitions},  Refname={Definition,Definitions}]  {Def}
\declaretheorem[style=scdef,sharenumber=Thm,     name=Remark,     refname={Remark,Remarks},          Refname={Remark,Remarks}]          {Rem}

\declaretheorem[style=scdef,sharenumber=Thm,     name=Example,    refname={Example,Examples},        Refname={Example,Examples}]        {Ex}
\declaretheorem[style=scprf,unnumbered,          name=Proof]{Prf}

\setlist[enumerate,1]{label={(\roman*)}}
\newcommand{\itemnr}[1]{(\romannumeral #1\relax)}

\DeclareTextFontCommand{\red}{\bfseries\color{red}}
\DeclareTextFontCommand{\green}{\bfseries\color{ForestGreen}}

\makeatletter\def\reallynopagebreak{\par\nopagebreak\@nobreaktrue}\makeatother

\def\svdots{\vbox{\baselineskip=4pt \lineskiplimit=0pt \kern2pt \hbox{.}\hbox{.}\hbox{.}\vspace{1pt}}}

\DeclareMathOperator{\Ad}{Ad}
\DeclareMathOperator{\ad}{ad}

\DeclareMathOperator{\pos}{pos}

\DeclareMathOperator{\SL}{SL}

\DeclareMathOperator{\SO}{SO}

\DeclareMathOperator{\U}{U}
\DeclareMathOperator{\SU}{SU}
\DeclareMathOperator{\Sp}{Sp}
\DeclareMathOperator{\PSp}{PSp}
\DeclareMathOperator{\Lag}{Lag}

\DeclareMathOperator{\Gr}{Gr}

\DeclareMathOperator{\Hom}{Hom}
\DeclareMathOperator{\Span}{span}

\DeclareMathOperator{\Sym}{Sym}
\DeclareMathOperator{\Is}{Is}
\DeclareMathOperator{\mbcd}{mbic}

\newcommand{\F}{\mathcal F}

\newcommand{\fg}{\mathfrak g}
\newcommand{\fk}{\mathfrak k}
\newcommand{\fp}{\mathfrak p}
\newcommand{\fa}{\mathfrak a}

\newcommand{\fn}{\mathfrak n}

\newcommand{\fso}{\mathfrak{so}}

\newcommand{\fsu}{\mathfrak{su}}
\newcommand{\aplusbar}{\overline{\fa^+}}

\newcommand{\RP}{\mathbb{R}\mathrm{P}}
\newcommand{\CP}{\mathbb{C}\mathrm{P}}

\newcommand{\bR}{\mathbb{R}}
\newcommand{\bC}{\mathbb{C}}
\newcommand{\bZ}{\mathbb{Z}}
\newcommand{\bN}{\mathbb{N}}
\newcommand{\bH}{\mathbb{H}}
\newcommand{\bK}{\mathbb{K}}

\newcommand{\bdry}{\partial_\infty\Gamma}
\newcommand{\ctheta}{\Delta\mathord\setminus\theta}
\newcommand{\ceta}{\Delta\mathord\setminus\eta}

\begin{document}

\title{Balanced ideals and domains of discontinuity of Anosov representations}
\author{Florian Stecker}

\maketitle

\begin{abstract}
  We consider the action of Anosov subgroups of a semi--simple Lie group on the associated flag manifolds. A systematic approach to construct cocompact domains of discontinuity for this action was given by Kapovich, Leeb and Porti in \cite{KapovichLeebPortiFlagManifolds}. For $\Delta$--Anosov representations, we prove that every cocompact domain of discontinuity arises from this construction, up to a few exceptions in low rank. Then we compute which flag manifolds admit these domains and, in some cases, the number of domains. We also find a new compactification for locally symmetric spaces arising from maximal representations into $\Sp(4n+2,\bR)$.
\end{abstract}


\section{Introduction}

Let $\Gamma$ be a word hyperbolic group and let $G$ be a semi--simple Lie group. A particularly well--behaved subset of the representations $\Hom(\Gamma,G)$ are the \emph{Anosov representations}. For instance, they have a discrete image and a finite kernel and they form an open subset of $\Hom(\Gamma, G)$. A definition can be found in \autoref{sec:prelim}. Examples of Anosov representations include all discrete injective representations into $\SL(2,\bR)$, quasi--Fuchsian representations into $\SL(2,\bC)$, representations in the Hitchin component of $\Hom(\pi_1S, \SL(n,\bR))$ for a closed surface $S$, and maximal representations from such a group $\pi_1S$ into a Hermitian Lie group.

We want to study the action of an Anosov representation $\rho$ on a flag manifold associated to $G$, that is a homogeneous space $\F = G/P$ where $P \subset G$ is a parabolic subgroup. A special case is the full flag manifold $\F_\Delta = G/B$, with $B$ being the minimal parabolic subgroup. If $G = \SL(n,\bR)$ the elements of flag manifolds are identified with sequences of nested subspaces of fixed dimensions in $\bR^n$. The $\rho$--action on $G/P$ is generally not proper, but in \cite{GuichardWienhardDomains} Guichard and Wienhard described a way of removing a ``bad set'' from a suitable flag manifold such that $\rho$ acts properly discontinuously and cocompactly on the complement. In other words, they constructed cocompact domains of discontinuity:

\begin{Def}
  A \emph{domain of discontinuity} $\Omega \subset \F$ for $\rho$ is a $\rho(\Gamma)$--invariant open subset such that the action $\Gamma \overset{\rho}{\curvearrowright} \Omega$ is proper. It is called \emph{cocompact} if the quotient $\Gamma \backslash \Omega$ is compact.
\end{Def}

Note that we require domains of discontinuity to be open subsets. In contrast, Danciger, Gueritaud and Kassel \cite{DancigerGueritaudKasselProjective,DancigerGueritaudKasselPseudoRiemannian} and Zimmer \cite{Zimmer} recently proved that the Anosov property is equivalent to the existence of certain cocompact domains in $\RP^n$ or $\bH^{p,q}$. These domains are closed subsets.

A systematic construction of (open) domains of discontinuity for Anosov representations was given by Kapovich, Leeb and Porti in \cite{KapovichLeebPortiFlagManifolds}. Say we have an Anosov representation $\rho$. It comes with a $\rho$--equivariant limit map $\xi \colon \bdry \to \F$ from the boundary of $\Gamma$ into some flag manifold $\F$. We want to find a cocompact domain of discontinuity in a flag manifold $\F'$, which may be different from $\F$. To construct such domains, Kapovich, Leeb and Porti use a combinatorial object called a \emph{balanced ideal}. That is a subset $I$ of the finite set $G \setminus (\F \times \F')$ satisfying certain conditions (see \autoref{sec:prelim}). They prove that, for every balanced ideal $I$, the set
\[\Omega_{\rho,I} = \F' \setminus \bigcup_{x\in\bdry} \{f \in \F' \mid G(\xi(x),f) \in I\}\]
is a cocompact domain of discontinuity.

In this paper, we prove a type of converse to this statement in the case of $\Delta$--Anosov representations. $\Delta$--Anosov, or Anosov with respect to the minimal parabolic, is the strongest form of the Anosov property, and far more restrictive than general Anosov representations.

\begin{Thm}\label{thm:intro_cocompact_domains_lift}
  Let $\rho \colon \Gamma \to G$ be a $\Delta$--Anosov representation and $\Omega \subset \F$ a cocompact domain of discontinuity for $\rho$ in some flag manifold $\F$. Then there is a balanced ideal $I \subset G \backslash (\F_\Delta \times \F_\Delta)$ such that the lift of $\Omega$ to the full flag manifold $\F_\Delta$ is a union of connected components of $\Omega_{\rho,I}$.
\end{Thm}

We can say more if the dimension of the bad set is not too big. To calculate it, we associate to a semi--simple Lie group $G$ a number $\mbcd(G)$ which gives a lower bound on the codimension of the set we have to remove for every limit point $x \in \bdry$. It increases with the rank of $G$, for example $\mbcd(\SL(n,\bR)) = \lfloor (n+1)/2\rfloor$. The general definition is given in \autoref{sec:dimensions}. With this we get

\begin{Thm}\label{thm:intro_cocompact_domains_correspondence}
  Let $\rho \colon \Gamma \to G$ be a $\Delta$--Anosov representation and assume that $\dim\bdry \leq \mbcd(G) - 2$. Then there is a 1:1 correspondence of balanced ideals in $G \backslash (\F_\Delta \times \F)$ and non--empty cocompact domains of discontinuity in $\F$.
\end{Thm}

\begin{Rem}
  A key point for these theorems is that cocompact domains are maximal among all domains of discontinuity, at least if they are connected. We then establish a correspondence between minimal fat ideals and maximal domains of discontinuity, even if they are not cocompact. This approach only works for $\Delta$--Anosov representations: \autoref{sec:counterexample} shows an example of an Anosov, but not $\Delta$--Anosov representation which admits infinitely many maximal domains of discontinuity. They are not cocompact.
\end{Rem}

\begin{Rem}
  For some choices of flag manifold $\F$ there are no balanced ideals in $G \backslash (\F_\Delta \times \F)$ and therefore no cocompact domain of discontinuity in $\F$. But even in these cases one can sometimes find a cocompact domain of discontinuity in a finite cover $\widehat\F$ of $\F$, an \emph{oriented flag manifold}. In \cite{SteckerTreib}, a theory analogous to \cite{KapovichLeebPortiFlagManifolds} is developed for oriented flag manifolds. Many arguments of this paper can potentially be extended to this setting, to show that essentially all cocompact domains in oriented flag manifolds are of the type described in \cite{SteckerTreib}.
\end{Rem}

Examples of $\Delta$--Anosov representations into $\SL(n,\bR)$ include

\begin{itemize}
\item representations in a \emph{Hitchin component} of $\Hom(\pi_1S,\SL(n,\bR))$ for a closed surface $S$. These are the connected components containing the composition of discrete injective representations $\pi_1S \to \SL(2,\bR)$ with the irreducible embedding $\SL(2,\bR) \to \SL(n,\bR)$.
\item the composition of a Fuchsian representation with the reducible embedding of $\SL(2,\bR)$ into $\SL(3,\bR)$ and small deformations thereof. These were studied in \cite{Barbot}. A similar construction works for all $\SL(n,\bR)$ with odd $n$.
\item representations of free groups $\rho \colon F_k \to \SL(2n,\bR)$ which arise as a sum of discrete injective representations $\rho_1,\dots,\rho_n \colon F_k \to \SL(2,\bR)$ such that $\rho_i$ uniformly dominates $\rho_j$ for all $i \geq j$. This means that $\log\lambda_1(\rho_i(\gamma)) \geq c\,\log\lambda_1(\rho_j(\gamma))$ for some constant $c > 1$ and all $\gamma \in F_k$, where $\lambda_1$ denotes the highest eigenvalue. Then $\rho$ is $\Delta$--Anosov, as are small deformations \cite[7.1]{GueritaudGuichardKasselWienhard}.
\end{itemize}

Next we combine \autoref{thm:intro_cocompact_domains_correspondence} with a criterion for the existence of balanced ideals in the case $G = \SL(n,\bR)$ or $G = \SL(n,\bC)$. This gives a full description which flag manifolds admit cocompact domains of discontinuity.

\begin{Thm}\label{thm:intro_Hitchin}
  Let $\Gamma$ be a hyperbolic group, $\bK \in \{\bR, \bC\}$, and $\rho \colon \Gamma \to \SL(n,\bK)$ a $\Delta$--Anosov representation. Choose integers $i_0, \dots,i_{k+1}$ with $0 = i_0 < i_1 < \dots < i_{k+1} = n$. Denote by $\F$ the corresponding flag manifold
  \[\F = \{V^{i_1} \subset \dots \subset V^{i_k} \subset \bK^n \mid \dim_\bK V^{i_j} = i_j\}.\]
  Assume that
  \[n \geq \begin{cases}2 \dim\bdry + 3 & \text{if $\bK = \bR$} \\ 2 \lfloor (\dim\bdry+1) / 2 \rfloor + 1 & \text{if $\bK = \bC$}\end{cases}\]
  and let $\delta = |\{0 \leq j \leq k \mid i_{j+1} - i_j\;\text{\upshape is odd}\}|$. Then
  \begin{enumerate}
  \item If $n$ is even, a non--empty cocompact open domain of discontinuity for $\Gamma \overset{\rho}{\curvearrowright} \F$ exists if and only if $\delta \geq 1$.
  \item If $n$ is odd, a non--empty cocompact open domain of discontinuity for $\Gamma \overset{\rho}{\curvearrowright} \F$ exists if and only if $\delta \geq 2$.
  \end{enumerate}
\end{Thm}

In particular, for surface group representations into $\SL(n,\bR)$ acting on Grassmannians we get

\begin{Cor}
  Let $n \geq 5$ and let $\rho \colon \pi_1S \to \SL(n, \bR)$ be a $\Delta$--Anosov representation from the fundamental group of a surface $S$ with or without boundary. Then the induced action $\pi_1S \overset{\rho}{\curvearrowright} \Gr(k,n)$ on the Grassmannian of $k$--planes in $\bR^n$ admits a non--empty cocompact domain of discontinuity if and only if $n$ is odd and $k$ is even. For small $n$, the number of different such domains is
  \begin{center}
  \begin{tabular}{l|ccccc}
    & $k=1$ & $k=3$ & $k=5$ & $k=7$ & $k=9$ \\
    \hline
    $n=6$ & $1$ & $2$ & $1$ && \\
    $n=8$ & $1$ & $7$ & $7$ & $1$ &\\
    $n=10$ & $1$ & $42$ & $2227$ & $42$ & $1$ \\
  \end{tabular}
\end{center}
\end{Cor}

\begin{Rem}
  While we don't know a general formula for these numbers, they are just the numbers of different balanced ideals. This is a combinatorial problem and needs no information about the representation except that it is $\Delta$--Anosov. So we can use a computer program to enumerate all balanced ideals. \autoref{sec:list_of_balanced_ideals} shows more results from this enumeration.
\end{Rem}

In low ranks, e.g. $\SL(n,\bR)$ with $n \leq 4$ if $\Gamma$ is a surface group, the existence of cocompact domains of discontinuity depends on more information about the geometry of $\rho$. So we can't make general lists like above in these cases. But at least for Hitchin representations the cocompact domains of discontinuity are also known in low ranks. We will briefly discuss this in \autoref{sec:Hitchin}.

Another observation from studying the list of balanced ideals is the existence of a balanced ideal in $\Sp(2n,\bC) \backslash (\Lag(\bC^{2n}) \times \Lag(\bC^{2n}))$ whenever $n$ is odd. This allows us to construct a compactification for locally symmetric spaces associated to $\{\alpha_n\}$--Anosov representations, which is modeled on the bounded symmetric domain compactification of the symmetric space. In particular, this includes maximal representations.

Recall that the symmetric space $X = \Sp(2n,\bR)/\U(n)$, like any Hermitian symmetric space, can be realized as a bounded symmetric domain $D \subset \bC^{n(n+1)/2}$ \cite[Theorem VIII.7.1]{Helgason}. Concretely, we can take as $D$ the set of symmetric complex matrices $Z$ such that $1 - \overline ZZ$ is positive definite. Its closure $\overline D$ in $\bC^{n(n+1)/2}$ is the \emph{bounded symmetric domain compactification} of $X$.

\begin{Thm}\label{thm:intro_compactification}
  Let $n$ be odd and $\rho \colon \Gamma \to \Sp(2n,\bR)$ be an $\{\alpha_n\}$--Anosov representation. Then there exists a subset $\widehat D \subset \bC^{n(n+1)/2}$ such that $D \subset \widehat D \subset \overline D$ on which $\rho(\Gamma)$ acts properly discontinuously with compact quotient. This quotient $\Gamma \backslash \widehat D$ is a compactification of the locally symmetric space $\Gamma \backslash D$.
\end{Thm}

\paragraph{Acknowledgements} I would like to thank David Dumas for developing a program to enumerate balanced ideals with me, Jean--Philippe Burelle for his help with the example in \autoref{sec:counterexample}, and Beatrice Pozzetti and Anna Wienhard for reading a draft of this paper. I am also grateful to Anna Wienhard and her whole research group in Heidelberg for countless helpful suggestions and for answering all kinds of questions.

{This project has received funding from the Klaus Tschira Foundation, the RTG 2229 grant of the German Research Foundation, the European Research Council (ERC) under the European Union's Horizon 2020 research and innovation programme (ERC consolidator grant 614733 (geometric structures), ERC starting grant 715982 (DiGGeS)), and the U.S. National Science Foundation grants DMS 1107452, 1107263, 1107367 ``RNMS: Geometric Structures and Representation Varieties'' (the GEAR Network). Part of the work was done while the author was visiting the Institut des Hautes Études Scientifiques (IHES).}

\section{Flag manifolds and Anosov representations}\label{sec:prelim}

In this section we fix some notation and conventions. In particular, we define Anosov representations and the relative position of flags.

Let $G$ be a connected semi--simple Lie group with finite center and $\fg$ its Lie algebra. Choose a maximal compact subgroup $K \subset G$ with Lie algebra $\fk \subset \fg$ and let $\fp = \fk^\perp$ be the Killing orthogonal complement in $\fg$. Further choose a maximal subspace $\fa \subset \fp$ on which the Lie bracket vanishes. Denote by $\fg_\alpha = \{X \in \fg \mid [H,X] = \alpha(H) X \ \forall H \in \fa\}$ the \emph{restricted root spaces} and $\Sigma = \{\alpha \in \fa^* \mid \fg_\alpha \neq 0\}$ the set of \emph{restricted roots}. Also choose a simple system $\Delta \subset \Sigma$ and let $\Sigma^\pm \subset \Sigma$ be the corresponding positive and negative roots. Note that any choice of the triple $(K,\fa,\Delta)$ is equivalent by conjugation in $G$ (see \cite[Theorem 2.1]{Helgason} and \cite[Theorems 2.63, 6.51, 6.57]{Knapp}).

The \emph{Weyl group} of $G$ is the group $W = N_K(\fa)/Z_K(\fa)$. It can be viewed as a group of linear isometries of $\fa$ equipped with the Killing form. A natural generating set of $W$ is given by $\Delta$, identifying every $\alpha \in \Delta$ with the orthogonal reflection along $\ker \alpha$. As $(W,\Delta)$ is a finite Coxeter system there is a unique \emph{longest element} $w_0 \in W$, which squares to the identity. The \emph{opposition involution} is the map $\iota(w) = w_0w w_0$ on $W$. It restricts to an involution $\iota \colon \Delta \to \Delta$ of the simple roots.

Define for $\varnothing \neq \theta \subset \Delta$ the Lie subalgebras
\[\fn = \bigoplus_{\alpha \in \Sigma^+} \fg_\alpha, \quad \fn^- = \bigoplus_{\alpha \in \Sigma^-} \fg_\alpha, \quad \fp_\theta = \bigoplus_{\!\!\!\!\alpha \in \Sigma^+ \cup \Span(\ctheta)\!\!\!\!} \fg_\alpha\]
and let $A, N, N^- \subset G$ be the connected Lie subgroups corresponding to $\fa$, $\fn$, and $\fn^-$. Note that $N^- = w_0 N w_0^{-1}$ (choosing any representative for $w_0$ in $N_K(\fa)$). Let $P_\theta = Z_G(\fp_\theta)$ be the \emph{parabolic subgroup} corresponding to $\theta$. The \emph{minimal parabolic subgroup} $B = P_\Delta$ decomposes as $B = Z_K(\fa) A N$ via the Iwasawa decomposition.

The \emph{flag manifolds} of $G$ are the spaces $\F_\theta = G/P_\theta$ for non--empty subsets $\theta \subset \Delta$. For two such subsets $\theta, \eta \subset \Delta$ the set of \emph{relative positions} is $W_{\theta,\eta} = \langle \ctheta \rangle \backslash W / \langle \ceta \rangle$, where $\langle A \rangle$ denotes the subgroup generated by $A \subset W$. The \emph{relative position map} is the unique map
\[\pos_{\theta,\eta} \colon \F_\theta \times \F_\eta \to W_{\theta,\eta}\]
which is invariant by the diagonal action of $G$ and satisfies $\pos_{\theta,\eta}([1],[w]) = w$ for all $w \in N_K(\fa)$. We will often omit the subscripts $\theta,\eta$. Two flags $f \in \F_\theta$ and $f' \in \F_\eta$ are said to be \emph{transverse} if $\pos_{\theta,\eta}(f,f') = w_0$.

For any $x \in \F_\theta$ and $w \in W_{\theta,\eta}$ the space
\[C_w(x) = \{y \in \F_\eta \mid \pos(x,y) = w\}\]
is the \emph{Bruhat cell} of $w$ relative to $x$ (it is indeed a cell if $\eta = \theta = \Delta$, but not in general). The \emph{Bruhat order} $\leq$ on $W_{\theta,\eta}$ is the inclusion order on the closures of Bruhat cells, i.e. for all $w, w' \in W_{\theta,\eta}$ and any (and therefore all) $x \in \F_\theta$
\[w \leq w' \quad \Leftrightarrow \quad \overline{C_w(x)} \subset \overline{C_{w'}(x)} \quad \Leftrightarrow \quad C_w(x) \subset \overline{C_{w'}(x)}.\]
Now assume that $\iota(\theta) = \theta$. Then $w_0$ acts on $W_{\theta,\eta}$ by left--multiplication.
\begin{Def}
  A subset $I \subset W_{\theta,\eta}$ is an \emph{ideal} if $w \in I$ implies $w' \in I$ whenever $w' \leq w$. It is \emph{fat} if $w \not\in I \Rightarrow w_0w \in I$ for all $w \in W_{\theta,\eta}$ and \emph{slim} if $w \in I \Rightarrow w_0w \not\in I$ for all $w \in W_{\theta,\eta}$. It is \emph{balanced} if it is fat and slim, i.e. $w \in I \Leftrightarrow w_0w \not\in I$.
\end{Def}

Let
\[\mu \colon G \to \aplusbar\]
be the \emph{Cartan projection} defined by the KAK--decomposition. That is, for every element $g \in G$ there are $k, \ell \in K$ and a unique $\mu(g) \in \aplusbar = \{X \in \fa \mid \alpha(X) \geq 0 \;\forall \alpha \in \Delta\}$ such that $g = k\,e^{\mu(g)}\,\ell$. $k$ and $\ell$ are uniquely defined up to an element in the centralizer of $\mu(g)$.

\begin{Def}\label{def:Anosov}
  Let $\theta \subset \Delta$ be non--empty and $\iota(\theta) = \theta$. A representation $\rho \colon \Gamma \to G$ is a \emph{$\theta$--Anosov representation} if there are constants $C,c > 0$ such that
  \[\alpha(\mu(\rho(\gamma))) \geq C|\gamma| - c \qquad \forall \alpha \in \theta, \gamma \in \Gamma,\]
  where $| \cdot |$ is the word length in $\Gamma$ with respect to any finite generating set.
\end{Def}

This is only one out of many equivalent definitions of Anosov representations \cite{Labourie,GuichardWienhardDomains,KapovichLeebPortiMorseActions,KapovichLeebPortiMorseLemma,GueritaudGuichardKasselWienhard,BochiPotrieSambarino}. See \cite[Theorem 1.1]{KapovichLeebPortiCharacterizations} for an overview and \cite{KapovichLeebPortiCharacterizations,BochiPotrieSambarino} for proofs of the equivalences.

The main fact we need about Anosov representations (which is part of many definitions) is that a $\theta$--Anosov representation $\rho \colon \Gamma \to G$ admits a unique \emph{limit map}
\[\xi \colon \bdry \to \F_\theta\]
which is continuous, $\rho$--equivariant, transverse (meaning that $\xi(x)$ and $\xi(y)$ are transverse for all $x \neq y$) and which maps the attracting fixed point of every infinite order element $\gamma \in \Gamma$ to the attracting fixed point of $\rho(\gamma)$.

\section{Domains of discontinuity}

In this section, we prove the main results, \autoref{thm:intro_cocompact_domains_lift} (which is \autoref{thm:cocompact_domains_lift}) and \autoref{thm:intro_cocompact_domains_correspondence} (as \autoref{thm:cocompact_domains_correspondence}), after some lemmas.

\subsection{Divergent sequences}

We first consider the behaviour of divergent sequences $(g_n) \in G^\bN$ in the semi--simple Lie group $G$. As before, let $\theta, \eta \subset \Delta$ be non--empty subsets of the simple restricted roots and assume $\iota(\theta) = \theta$.

\begin{Def}
  A sequence $(g_n) \in G^\bN$ is \emph{$\theta$--divergent} if $\alpha(\mu(g_n)) \to \infty$ for all $\alpha \in \theta$.
\end{Def}

\begin{Def}
  A $\theta$--divergent sequence $(g_n) \in G^\bN$ is \emph{simply $\theta$--divergent} if it has KAK--de\-compo\-si\-tions $g_n = k_n a_n \ell_n$ such that $(k_n)$ and $(\ell_n)$ converge to some $k,\ell \in K$. Then $g^- = [\ell^{-1}w_0] \in \F_\theta$ and $g^+ = [k] \in \F_\theta$ are the \emph{repelling and attracting limits} of this sequence.
\end{Def}

\begin{Rem}
  By compactness of $K$, every $\theta$--divergent sequence in $G$ has a simply $\theta$--divergent subsequence. The limits $(g^-,g^+)$ do not depend on the choice of decomposition. If a simply $\theta$--divergent sequence $(g_n) \in G^\bN$ has limits $(g^-, g^+) \in \F_\theta^2$, then $(g_n^{-1})$ is also simply $\theta$--divergent with limits $(g^+, g^-)$.
\end{Rem}

In \cite{KapovichLeebPortiFlagManifolds}, the following characterization of $\theta$--divergent sequences is used.

\begin{Lem}\label{lem:theta_divergent_implies_flag_convergence}
  Let $(g_n) \in G^\bN$ be a simply $\theta$--divergent sequence with limits $(g^-,g^+) \in \F_\theta^2$. Then
  \[g_n|_{C_{w_0}(g^-)} \to g^+\]
  locally uniformly as functions from $\F_\theta$ to $\F_\theta$ (where $g^+$ is the constant function).
\end{Lem}

\begin{Prf}
  By assumption, we can write $g_n = k_n e^{A_n} \ell_n$ with $k_n \to k$ and $\ell_n \to \ell$ in $K$, $A_n \in \aplusbar$ and $\alpha(A_n) \to \infty$ for all $\alpha \in \theta$. Furthermore, $g^- = [\ell^{-1}w_0]$ and $g^+ = [k]$. Now let $(f_n) \in \F_\theta^\bN$ be a sequence converging to $f \in C_{w_0}(g^-)$. Then $\ell_n f_n \to \ell f \in C_{w_0}([w_0])$. Since this is an open set, we can assume that $\ell_n f_n \in C_{w_0}([w_0])$ for all $n$. By the Langlands decomposition of $P_\theta$ \cite[Proposition 7.83]{Knapp} we can write $\ell_n f_n = [\exp(X_n)]$ with
  \[X_n = \sum_{\alpha} X_n^\alpha \in \bigoplus_{\!\!\!\alpha \in \Sigma^- \setminus \Span(\ctheta)\!\!\!} \fg_\alpha.\]
  All of the roots $\alpha$ appearing in this sum are linear combinations of simple roots with only non--positive coefficients, and with at least one coefficient of a root in $\theta$ being negative. So $\alpha(A_n) \to -\infty$ for all such roots $\alpha$, and therefore
  \[g_n f_n = [k_n e^{A_n} e^{X_n}] = [k_n e^{A_n} e^{X_n} e^{-A_n}] = \left[k_n \exp\left(\sum_\alpha e^{\alpha(A_n)}X_n^\alpha\right)\right] \to [k] = g^+. \qedhere\]
\end{Prf}

For the action of a discrete group $\Gamma$ on a manifold $X$, there is a useful reformulation of properness. By \cite[Proposition 1]{Frances}, the action is proper if and only if it has no dynamical relations, in the following sense:

\begin{Def}
  Let $\Gamma$ be a discrete group acting smoothly on a manifold $X$. Let $(\gamma_n) \in \Gamma^\bN$ be a divergent sequence (i.e. no element occurs infinitely many times) and $x,y \in X$. Then \emph{$x$ is dynamically related to $y$ via $(\gamma_n)$}, $x \stackrel{(\gamma_n)}{\sim} y$, if there is a sequence $(x_n) \in X^\bN$ such that
  \[x_n \to x\quad \text{and} \quad \gamma_n x_n \to y.\]
  We say $x$ and $y$ are \emph{dynamically related}, $x \sim y$, if they are dynamically related via any divergent sequence in $\Gamma$.
\end{Def}

\Autoref{lem:relative_position_inequality} is the key step of the proof of proper discontinuity in \cite{KapovichLeebPortiFlagManifolds} (Proposition 6.5). It states that flags in $\F_\eta$ can only be dynamically related by the action of $\rho$ if their relative positions satisfy the inequality \eqref{eq:dynamical_relation_inequality}. \Autoref{lem:deltanecessary} is a converse to this statement in the case $\theta = \Delta$: It says that whenever two flags $f,f'$ satisfy a relation like \eqref{eq:dynamical_relation_inequality}, then they are indeed dynamically related.

\begin{Lem}\label{lem:relative_position_inequality}
  Let $(g_n) \in G^\bN$ be a simply $\theta$--divergent sequence and let $(g^-,g^+) \in \F_\theta^2$ be its limits. Let $f, f' \in \F_\eta$ be dynamically related via $(g_n)$. Then
  \begin{equation}
    \label{eq:dynamical_relation_inequality}
    \pos_{\theta,\eta}(g^+, f') \leq w_0 \pos_{\theta,\eta}(g^-, f).
  \end{equation}
\end{Lem}

\begin{Prf}
  As $f, f'$ are dynamically related, there exists a sequence $f_n \in \F_\eta$ converging to $f$ such that $g_nf_n \to f'$. We can write $f_n = h_n f$ for some sequence $h_n \in G$ converging to $1$. Let $w = \pos_{\theta,\eta}(g^-, f)$. Then there exists $g \in G$ such that $g(g^-, f) = ([1], [w])$. Define $F = [g^{-1} w_0] \in \F_\theta$. We get the following relative positions:
  \[\pos_{\theta, \theta}(g^-, F) = w_0, \qquad \pos_{\theta, \eta}(g^-, f) = w, \qquad \pos_{\theta, \eta}(F, f) = w_0w.\]
  Now $h_n F \to F$ and since $F$ and almost all of the $h_nF$ are in $C_{w_0}(g^-)$, we get by \autoref{lem:theta_divergent_implies_flag_convergence} that $g_nh_nF \to g^+$. So
  \[\pos(g^+, f') \leq \pos(g_n h_n F, g_n h_n f) = \pos(F, f) = w_0 w = w_0 \pos(g^-, f).\qedhere\]
\end{Prf}

\begin{Lem}\label{lem:unipotent_limits}
  Let $(A_n) \in \overline{\mathfrak a^+}^\bN$ be a $\Delta$--divergent sequence and $n^+ \in N, n^- \in N^-$. Then there exists a sequence $(h_n) \in G^\bN$ such that $h_n \to n^-$ and $e^{A_n}h_ne^{-A_n} \to n^+$.
\end{Lem}

\begin{Prf}
  We can write $n^- = e^{X^-}$ for $X^- \in \mathfrak n^-$ and $n^+ = e^{X^+}$ with $X^+ \in \mathfrak n$. Let $H_n = X^- + e^{-\ad A_n} X^+$ and $h_n = e^{H_n}$. For all $\alpha \in \Sigma^+$ and $X_\alpha \in \mathfrak g_\alpha$ we know that $e^{-\ad A_n} X_\alpha = e^{-\alpha(A_n)} X_\alpha$ converges to $0$. As $X^+$ is a linear combination of these, $H_n \to X^-$ and thus $h_n \to n^-$. On the other hand
  \[e^{A_n} e^{H_n} e^{-A_n} = \exp(\Ad_{e^{A_n}} H_n) = \exp(e^{\ad A_n}H_n) = \exp(e^{\ad A_n}X^- + X^+)\]
  which converges to $n^+ = e^{X^+}$ by a similar argument.
\end{Prf}

\begin{Lem}\label{lem:deltanecessary}
  Let $(g_n) \in G^\bN$ be a simply $\Delta$--divergent sequence with limits $(g^-,g^+) \in \F_\Delta^2$. Let $f,f' \in \F_\eta$ and $w \in W_{\Delta,\eta}$ with
  \begin{equation}
    \label{eq:relative_position_equality}
    \pos_{\Delta,\eta}(g^-,f) = w, \qquad \pos_{\Delta,\eta}(g^+,f') = w_0w.
  \end{equation}
  Then $f$ is dynamically related to $f'$ via $(g_n)$.
\end{Lem}

\begin{Prf}
  Fix some representative in $N_K(\mathfrak a)$ for $w$ and $w_0$. Let $g_n = k_n e^{A_n} \ell_n$ be a KAK--decomposition, such that $(k_n), (\ell_n) \in K^\bN$ converge to $k,\ell \in K$ and $\alpha(A_n) \to \infty$ for all $\alpha \in \Delta$. Then the limits can be written as
  \[g^- = [\ell^{-1}w_0] \in \F_\Delta, \qquad g^+ = [k] \in \F_\Delta.\]
  Because of \eqref{eq:relative_position_equality} there exist $h,h' \in G$ with
  \[g^- = [h] \in \F_\Delta, \qquad f = [hw] \in \F_\eta, \qquad g^+ = [h'] \in \F_\Delta, \qquad f' = [h'w_0w] \in \F_\eta.\]
  So $w_0^{-1} \ell h, k^{-1}h' \in B$, which means we can write $w_0^{-1} \ell h = n a m$ and $k^{-1} h' = n' a' m'$ for some $n,n' \in N$, $a,a' \in A$ and $m,m' \in Z_K(\mathfrak a)$. Consequently,
  \begin{align*}
    f &= [hw] = [\ell^{-1}w_0namw] = [\ell^{-1} (w_0 n w_0^{-1}) w_0 w], \\
    f' &= [h'w_0w] = [kn'a'm'w_0w] = [kn'w_0w],
  \end{align*}
  since elements of $A$ and $Z_K(\mathfrak a)$ commute with Weyl group elements. By \autoref{lem:unipotent_limits} there is a sequence $(h_n) \in G^\bN$ such that $h_n \to w_0 n w_0^{-1}$ and $e^{A_n}h_n e^{-A_n} \to n'$. Let $f_n = [\ell_n^{-1}h_nw_0w] \in \F_\eta$. Then $f_n \to f$ and
  \[g_nf_n = [k_ne^{A_n} h_n w_0 w] = [k_n e^{A_n} h_n e^{-A_n} w_0 w] \to [k n' w_0 w] = f'. \qedhere\]
\end{Prf}

\subsection{Limit sets}

Now let $\Gamma$ be a non--elementary hyperbolic group and $\rho \colon \Gamma \to G$ a representation.

\begin{Def}
  Let $\theta \subset \Delta$ be non--empty and $\iota$--invariant. The \emph{limit set} of $\rho$ is the set
  \[\Lambda_{\rho,\theta} = \{g^+ \mid \exists \text{$(g_n) \in \rho(\Gamma)^\bN$ simply $\theta$--divergent with limits $(g^-,g^+)$}\} \subset \F_\theta\]
  and the set of \emph{limit pairs} is
  \[\Lambda_{\rho,\theta}^{[2]} = \{(g^-,g^+) \mid \text{$\exists (g_n) \in \rho(\Gamma)^\bN$ simply $\theta$--divergent with limits $(g^-,g^+)$}\} \subset \F_\theta^2.\]
\end{Def}

These limit sets are particularly well--behaved for Anosov representations. Namely, we have the following well--known facts:

\begin{Prop}\label{pro:lambda2xi}
  If $\rho$ is $\theta$--Anosov with limit map $\xi \colon \bdry \to \F_\theta$, then
  \[\Lambda_{\rho,\theta} = \xi(\bdry), \qquad \Lambda_{\rho,\theta}^{[2]} = \xi(\bdry)^2.\]
\end{Prop}

The first part can be found e.g. in \cite[Theorem 5.3(3)]{GueritaudGuichardKasselWienhard}. We give a detailed proof of the second part. We first need two short lemmas.

\begin{Lem}\label{lem:boundary_convergence_implies_flag_convergence}
  Let $\rho \colon \Gamma \to G$ be a $\theta$--Anosov representation with limit map $\xi \colon \bdry \to \F_\theta$, $(\gamma_n) \in \Gamma^\bN$ a diverging sequence and $(\gamma^-,\gamma^+) \in \bdry^2$ such that
  \[\gamma_n|_{\bdry \setminus \{\gamma^-\}} \to \gamma^+\]
  locally uniformly. Then on $\F_\theta$ we also have the locally uniform convergence
  \begin{equation}
    \label{eq:limits_in_F}
    \rho(\gamma_n)|_{C_{w_0}(\xi(\gamma^-))} \to \xi(\gamma^+).
  \end{equation}
\end{Lem}

\begin{Prf}
  By restricting to a subsequence we can assume that $\rho(\gamma_n)$ is simply $\theta$--divergent with limits $(g^-,g^+) \in \F_\theta^2$, so $\rho(\gamma_n)|_{C_{w_0}(g^-)} \to g^+$ by \autoref{lem:theta_divergent_implies_flag_convergence}. Since $\Lambda_{\rho,\theta} = \xi(\bdry)$, we have $g^- = \xi(x)$ for some $x \in \bdry$. Let $z \in \bdry \setminus \{\gamma^-, x\}$. Then $\xi(z) \in C_{w_0}(\xi(x))$, so $\rho(\gamma_n z) = \rho(\gamma_n)\xi(z) \to g^+$, but also $\gamma_n z \to \gamma^+$, so $g^+ = \xi(\gamma^+)$. The same argument applied to $\gamma_n^{-1}$ instead of $\gamma_n$ shows that $g^- = \xi(\gamma^-)$. Finally, since any subsequence of $(\gamma_n)$ has a subsequence satisfying \eqref{eq:limits_in_F}, it actually holds for the whole sequence $(\gamma_n)$.
\end{Prf}

\begin{Lem}\label{lem:limit_pairs_in_Gamma}
  Let $(\gamma^-,\gamma^+) \in \bdry^2$. Then there exists a divergent sequence $(\gamma_n) \in \Gamma^\bN$ such that
  \begin{equation}
    \label{eq:limits_in_Gamma}
    \gamma_n|_{\bdry \setminus \{\gamma^-\}} \to \gamma^+
  \end{equation}
  locally uniformly.
\end{Lem}

\begin{Prf}
  Fix a metric on $\bdry$. Let $P \subset \bdry^2$ be the set of fixed point pairs of infinite order elements of $\Gamma$. $P$ is dense in $\bdry^2$ by \cite[8.2.G]{Gromov}, so we find a sequence $(\gamma_n) \in \Gamma^\bN$ whose fixed point pairs $(\gamma_n^-,\gamma_n^+)$ approach $(\gamma^-,\gamma^+)$. Substituting each $\gamma_n$ by a sufficiently high power, we can assume that $\gamma_n$ maps the complement of $B_{1/n}(\gamma_n^-)$ into $B_{1/n}(\gamma_n^+)$ by \cite[8.1.G]{Gromov}. Now let $x_n \to x$ be any convergent sequence in $\bdry$ with $x \neq \gamma^-$. Once $n$ is large enough such that
  \[\textstyle d(x_n,x) \leq \frac{1}{3}d(x,\gamma^-), \quad d(\gamma_n^-, \gamma^-) \leq \frac{1}{3}d(x,\gamma^-), \quad 1/n \leq \frac{1}{3}d(x,\gamma^-),\]
  then $x_n$ lies outside of $B_{1/n}(\gamma_n^-)$, so $\gamma_nx_n \in B_{1/n}(\gamma_n^+)$ and therefore $\gamma_nx_n \to \gamma^+$.
\end{Prf}

\begin{Prf}[of $\Lambda_{\rho,\theta}^{[2]} = {\xi(\bdry)}^2$]
  We first prove $\xi(\bdry)^2 \subset \Lambda_{\rho,\theta}^{[2]}$. Let $(\gamma^-,\gamma^+) \in \bdry^2$ and let $(\gamma_n) \in \Gamma^\bN$ be a divergent sequence satisfying \eqref{eq:limits_in_Gamma}, which exists by \autoref{lem:limit_pairs_in_Gamma}. By \autoref{lem:boundary_convergence_implies_flag_convergence} this implies
  \[\rho(\gamma_n)|_{C_{w_0}(\xi(\gamma^-))} \to \xi(\gamma^+).\]
  Let $\rho(\gamma_{n_k})$ be any simple subsequence of the $\theta$--divergent sequence $(\rho(\gamma_n))$ and $(g^-,g^+) \in \F_\theta^2$ its limits. Then by \autoref{lem:theta_divergent_implies_flag_convergence}
  \[\rho(\gamma_{n_k})|_{C_{w_0}(g^-)} \to g^+,\]
  so $g^+ = \xi(\gamma^+)$ since $C_{w_0}(\xi(\gamma^-)) \cap C_{w_0}(g^-) \neq \varnothing$. It is not hard to see that \eqref{eq:limits_in_Gamma} implies $\gamma_n^{-1}|_{\bdry \setminus \gamma^+} \to \gamma^-$, and repeating the argument for this sequence shows that $g^- = \xi(\gamma^-)$. So $(\xi(\gamma^-),\xi(\gamma^+)) \in \Lambda_{\rho,\theta}^{[2]}$.

  For the other direction, let $(g^-,g^+) \in \Lambda_{\rho,\theta}^{[2]}$. Then there exists a divergent sequence $(\gamma_n) \in \Gamma^\bN$ such that $(\rho(\gamma_n))$ is simply $\theta$--divergent with limits $(g^-,g^+)$. Passing to a subsequence, we can also assume that it satisfies \eqref{eq:limits_in_Gamma} for some pair $(\gamma^-,\gamma^+) \in \bdry^2$, since $\Gamma \curvearrowright \bdry$ is a convergence group action. As before, we can use \autoref{lem:boundary_convergence_implies_flag_convergence} as well as \autoref{lem:theta_divergent_implies_flag_convergence} to show that $g^+ = \xi(\gamma^+)$, and $g^- = \xi(\gamma^-)$ by applying the same argument to the sequence of inverses. So $(g^-,g^+) \in \xi(\bdry)^2$.
\end{Prf}

\subsection{Maximal domains of discontinuity}

Recall that we call $\Omega \subset \F_\eta$ a domain of discontinuity if it is an open $\Gamma$--invariant subset on which $\Gamma$ acts properly. In this section, we deal with maximal domains of discontinuity, i.e. those which are not contained in any strictly larger domain of discontinuity.

\begin{Def}
  Let $\Omega \subset \F_\eta$,  $\Lambda \subset \F_\theta$, and $I \subset W_{\theta,\eta}$. We define
  \begin{align*}
  \Omega(\Lambda,I) &= \{x \in \F_\eta \mid \pos(\ell,x) \not\in I \ \forall \ell \in \Lambda\} \subset \F_\eta \\
  I(\Lambda,\Omega) &= W_{\theta,\eta} \setminus \{\pos(\ell,x) \mid \ell \in \Lambda, x \in \Omega \} \subset W_{\theta,\eta}
  \end{align*}
\end{Def}

\begin{Rem}
  It is easy to see that $I \subset I(\Lambda,\Omega(\Lambda, I))$ and $\Omega \subset \Omega(\Lambda, I(\Lambda, \Omega))$. If $I$ is an ideal and $\Lambda$ is closed, then $\Omega(\Lambda, I)$ is open (see e.g. \cite[Lemma 3.35 \& Lemma B.8(ii)]{SteckerTreib} for a detailed proof of this topological fact).
\end{Rem}

\begin{Prop}\label{pro:maximal_domains}
  Let $\rho \colon \Gamma \to G$ be $\Delta$--Anosov with limit map $\xi \colon \bdry \to \F_\Delta$ and $\Lambda = \xi(\bdry)$. Let $\Omega \subset \F_\eta$ be a maximal domain of discontinuity of $\rho$. Then $I \coloneqq I(\Lambda,\Omega) \subset W_{\Delta,\eta}$ is a fat ideal and $\Omega = \Omega(\Lambda, I)$.
\end{Prop}

\begin{Prf}
  We first prove that $I$ is an ideal. If not, there are $w' \leq w$ with $w \in I$ and $w' \not\in I$. So there exist $\ell \in \Lambda$ and $x \in \Omega$ such that $\pos(\ell,x) = w'$, i.e. $x \in C_{w'}(\ell)$. But $C_w(\ell) \subset \F_\eta \setminus \Omega$ which is closed, so $x \in C_{w'}(\ell) \subset \overline{C_w(\ell)} \subset \F_\eta \setminus \Omega$, a contradiction. So $I$ is an ideal.

  If $I$ was not fat, there would be a $w \in W_{\Delta,\eta}$ with $\pos(\ell,x) = w$ and $\pos(\ell',x') = w_0 w$ for some $\ell,\ell' \in \Lambda$ and $x,x' \in \Omega$. But $(\ell,\ell') \in \Lambda_{\rho,\Delta}^{[2]}$ by \autoref{pro:lambda2xi}, so there is a simply $\Delta$--divergent sequence $(\rho(\gamma_n)) \in \rho(\Gamma)^\bN$ with limits $(\ell,\ell')$. So $x$ and $x'$ would be dynamically related by \autoref{lem:deltanecessary}, a contradiction.

  Now $\Omega \subset \Omega(\Lambda, I)$ and the main theorem of \cite{KapovichLeebPortiFlagManifolds} says that $I$ being fat implies $\Omega(\Lambda, I)$ is a domain of discontinuity. So by maximality $\Omega = \Omega(\Lambda, I)$.
\end{Prf}


\begin{Cor}\label{cor:maximal_domains2}
  Every maximal domain of discontinuity of a $\Delta$--Anosov representation with limit set $\Lambda$ is of the form $\Omega(\Lambda, I)$ for a minimal fat ideal $I \subset W_{\Delta,\eta}$.
\end{Cor}

\begin{Prf}
  Let $\Omega$ be a maximal domain of discontinuity. Then by \autoref{pro:maximal_domains} $I(\Lambda, \Omega)$ is a fat ideal and $\Omega = \Omega(\Lambda, I(\Lambda, \Omega))$. In general, there could be other ideals generating the same domain. Let $I$ be minimal among all fat ideals $\widetilde I$ with $\Omega(\Lambda, \widetilde I) = \Omega$. Then $I$ is in fact minimal among all fat ideals, as otherwise there would be another fat ideal $I'$ with $\Omega = \Omega(\Lambda,I) \subsetneq \Omega(\Lambda,I')$, contradicting maximality of $\Omega$.
\end{Prf}

\subsection{Cocompactness}

The most important fact we need about cocompact domains of discontinuity is that they are essentially maximal. More precisely:

\begin{Lem}\label{lem:cocompact_maximal}
  Let $\rho \colon \Gamma \to G$ be a representation and $\Omega \subset \F_\eta$ a cocompact domain of discontinuity. Then $\Omega$ is a union of connected components of a maximal domain of discontinuity.
\end{Lem}

\begin{Prf}
  By Zorn's lemma $\Omega$ is contained in some maximal domain of discontinuity $\widetilde\Omega \in \F_\eta$ and it is an open subset. Then also $\Gamma \backslash \Omega \subset \Gamma \backslash \widetilde\Omega$, where $\Gamma \backslash \Omega$ is compact and $\Gamma \backslash \widetilde \Omega$ is Hausdorff. So $\Gamma \backslash \Omega$ is closed in $\Gamma \backslash \widetilde \Omega$ and therefore $\Omega$ is also closed in $\widetilde \Omega$.
\end{Prf}

This immediately leads to our first main theorem.

\begin{Thm}\label{thm:cocompact_domains_lift}
  Let $\rho \colon \Gamma \to G$ be a $\Delta$--Anosov representation with limit map $\xi \colon \bdry \to \F_\Delta$ and let $\Omega \subset \F_\eta$ be a cocompact domain of discontinuity for $\rho$. Then there is a balanced ideal $\widetilde I \subset W$ such that $\pi_\eta^{-1}(\Omega)$ is a $\Gamma$--invariant union of connected components of $\Omega(\xi(\bdry),\widetilde I) \subset \F_\Delta$.
\end{Thm}

\begin{Prf}
  The natural projection $\pi_\eta \colon \F_\Delta \to \F_\eta$ is smooth, $G$--equivariant and proper. This implies that $\widetilde\Omega = \pi_\eta^{-1}(\Omega)$ is also a cocompact domain of discontinuity. So by \autoref{lem:cocompact_maximal} there is a maximal domain of discontinuity $\widehat \Omega \subset \F_\Delta$ and $\widetilde \Omega$ is a union of connected components of $\widehat \Omega$. By \autoref{cor:maximal_domains2} $\widehat\Omega = \Omega(\Lambda, I)$ for a minimal fat ideal $I \subset W$. But since the action of $w_0$ on $W$ has no fixed points, every minimal fat ideal in $W$ is balanced. This is proved e.g. in \cite[Lemma 3.34]{SteckerTreib}.
\end{Prf}

We know from \cite{KapovichLeebPortiFlagManifolds} that domains constructed from a balanced ideal are cocompact. The combination of the next two lemmas shows that if the domain is dense, the converse also holds. That is, if a domain constructed from a fat ideal is cocompact, then this ideal must be balanced.

\begin{Lem}\label{lem:not_cocompact}
  Let $I \subset W_{\theta,\eta}$ be a fat ideal and $\Lambda$ the limit set of a $\theta$--Anosov representation $\rho$. Let
  \[D(\Lambda, I) = \{x \in \F_\eta \mid \exists \ell \neq \ell' \in \Lambda \colon \pos(\ell,x), \pos(\ell',x) \in I\}\]
  and let $\Omega_0 \subset \Omega$ be a $\rho(\Gamma)$--invariant union of connected components of $\Omega \coloneqq \Omega(\Lambda, I)$. Then $\Omega_0$ can be cocompact only if \/ $\overline{\Omega_0} \cap D(\Lambda, I) = \varnothing$.
\end{Lem}

\begin{Prf}
  Assume that $\Omega_0$ is cocompact and $x \in \partial\Omega_0$. Take a sequence $(x_n) \in \Omega_0^\bN$ with $x_n \to x$. Let $(h_n) \in G^\bN$ be a sequence converging to the identity such that $x_n = h_n x$. By cocompactness, a subsequence of $(x_n)$ converges in the quotient. Passing to this subsequence, there is $(g_n) \in \rho(\Gamma)^\bN$ such that $g_nx_n \to x' \in \Omega_0$. Clearly $g_n \to \infty$ as otherwise a subsequence of $(g_nx_n)$ would converge to something in $\partial\Omega_0$. Passing to a subsequence another time we can also assume that $(g_n)$ is simply $\theta$--divergent with limits $(g^-,g^+) \in \Lambda^2$.

  Now let $\ell \in \Lambda \setminus \{g^-\}$. Then $h_n \ell \to \ell$ and thus $g_n h_n \ell \to g^+$ by \autoref{lem:theta_divergent_implies_flag_convergence} since $\ell \in C_{w_0}(g^-)$ and this is an open set. So
  \[\pos(g^+,x') \leq \pos(g_n h_n \ell, g_n h_n x) = \pos(\ell, x).\]
  Since $x' \in \Omega(\Lambda, I)$ we know $\pos(g^+,x') \not\in I$, so $\pos(\ell, x) \not\in I$. This holds for every $\ell \in \Lambda \setminus \{g^-\}$, so $x \not\in D(\Lambda, I)$. We have thus proved that $\partial\Omega_0 \cap D(\Lambda,I) = \varnothing$. Also $\Omega_0 \cap D(\Lambda, I) = \varnothing$ holds by definition.
\end{Prf}

\begin{Lem}\label{lem:slim_no_intersections}
  In the setting of \autoref{lem:not_cocompact} an ideal $I \subset W_{\theta,\eta}$ is slim if and only if $D(\Lambda, I) = \varnothing$.
\end{Lem}

\begin{Prf}
  First assume that $I$ is not slim, i.e. there is $w \in I$ with $w_0 w \in I$. Let $\ell \neq \ell' \in \Lambda$. Since $\ell,\ell'$ are transverse there is $g \in G$ such that $g\ell = [1]$ and $g\ell' = [w_0] \in \F_\theta$. Let $x = [g^{-1} w] \in \F_\eta$. Then
  \[\pos_{\theta,\eta}(\ell, x) = \pos_{\theta,\eta}([1], [w]) = w \in I, \quad \pos_{\theta,\eta}(\ell', x) = \pos_{\theta, \eta}([w_0], [w]) = w_0w \in I,\]
  so $x \in D(\Lambda, I)$. Conversely, suppose that $x \in D(\Lambda, I)$. Then there are transverse $\ell, \ell' \in \Lambda \subset \F_\theta$ such that $\pos(\ell, x), \pos(\ell', x) \in I$. Let $g \in G$ with $g\ell = [1]$ and $g\ell' = [w_0]$. As for any flag, there exist $n \in N$ and $w \in W$ with $gx = [nw] \in \F_\eta$. Choose any $\Delta$--divergent sequence $a_k = e^{A_k} \in A$ with $A_k \in \overline{\mathfrak a^+}$. Then $a_k^{-1} n a_k \to 1$. Now
  \[\pos_{\theta,\eta}([w_0],[nw]) = \pos_{\theta,\eta}([w_0 w_0^{-1} a_k w_0],[nw w^{-1} a_k w]) = \pos_{\theta,\eta}([w_0 ],[a_k^{-1} n a_k w])\]
  and since $[a_k^{-1} n a_k w] \to [w]$ we get $\pos(\ell',x) = \pos([w_0],[nw]) \geq w_0 w$ and thus $w_0 w \in I$. But also $w = \pos([1],[nw]) = \pos(\ell, x) \in I$, so $I$ is not slim.
\end{Prf}

\subsection{Dimensions}\label{sec:dimensions}

If the domain $\Omega$ comes from a balanced ideal, the ``bad set'' $\F_\eta \setminus \Omega$ fibers over $\bdry$. The dimension of the fiber is bounded by the following quantity, depending only on $G$:

\begin{Def}
  For a subset $A \subset \Sigma$ of the simple roots let
  \[\dim A = \sum_{\alpha \in A} \dim \fg_\alpha.\]
  The we can define the \emph{minimal balanced ideal codimension} of $G$
  \[\mbcd(G) = \min_{w \in W \atop w_0 w \not\geq w} \dim\Psi_w, \qquad \Psi_w = \Sigma^+ \cap w \Sigma^-.\]
\end{Def}

Dumas and Sanders showed in \cite[Theorem 4.1]{DumasSanders} that if the Weyl group $W$ of $G$ has no factors of type $A_1$, then $w_0 w \leq w$ for all $w \in W$ with $\ell(w) \leq 1$, and that the same is true for $\ell(w) \leq 2$ if $W$ also has no factors of type $A_2$, $A_3$ or $B_2$. This implies $\mbcd(G) \geq 2$ resp. $\mbcd(G) \geq 3$ in these cases (and even higher lower bounds if the root spaces are more than one--dimensional, e.g. in the case of complex groups).

\begin{Ex}\label{ex:mbcd_SLn}
  For the special linear group we have
  \[\mbcd(\SL(n,\bR)) = \left\lfloor \frac{n+1}{2} \right\rfloor, \qquad \mbcd(\SL(n,\bC)) = 2\left\lfloor \frac{n+1}{2} \right\rfloor.\]
  To see this, recall that the Weyl group of $\SL(n,\bR)$ can be identified with the symmetric group $S_n$ with its standard generating set of adjacent transpositions. There is also a simple description of the Bruhat order on $S_n$: Define, for any permutation $w \in S_n$ and integers $i,j$
  \[w[i,j] \coloneqq |\{a \leq i \mid w(a) \leq j\}|.\]
  Then $w \leq w'$ if and only if $w[i,j] \geq w'[i,j]$ for all $i,j$ \cite[Theorem 2.1.5]{BjoernerBrenti}. So $w \leq w_0w$ if and only if
  \[w[i,j] \geq w_0w[i,j] = |\{a \leq i \mid w(a) > n-j\}|\]
  for all $i,j$. Since every root space $\fg_\alpha$ is 1--dimensional and $|\Psi_w| = \ell(w)$, $\mbcd(G)$ is therefore the minimal word length an element $w \in S_n$ has to have such that there are $i,j$ not satisfying this inequality. We can express this problem in a nice graphical way: Suppose we have $n$ balls in a row which we can permute. What is the minimal length of a permutation such that for some choice of $i$ and $j$, if the first $i$ balls were painted red before, then after the permutation there are more red ones among the last $j$ than among the first $j$?
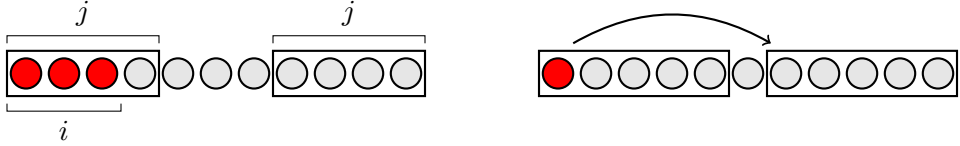
  \begin{center}
    \begin{tikzpicture}
      \foreach \i in {0,1,2} {
        \draw[thick,fill=red] (\i*0.5,0) circle (0.2);
      };
      \foreach \i in {3,4,5,6,7,8,9,10} \draw[thick,fill=black!10] (\i*0.5,0) circle (0.2);
      \draw[thick] (-0.25,-0.3) rectangle (-0.25 + 0.5*4,0.3);
      \draw[thick] (-0.25 + 0.5*7,-0.3) rectangle (-0.25 + 0.5*11,0.3);
      \draw (-0.25,0.4) -- (-0.25,0.5) -- (-0.25 + 0.5*4,0.5) node[midway,above] {$j$} -- (-0.25 + 0.5*4,0.4);
      \draw (-0.25 + 0.5*7,0.4) -- (-0.25 + 0.5*7,0.5) -- (-0.25 + 0.5*11,0.5) node[midway,above] {$j$} -- (-0.25 + 0.5*11,0.4);
      \draw (-0.25,-0.4) -- (-0.25,-0.5) -- (-0.25 + 0.5*3,-0.5) node[midway,below] {$i$} -- (-0.25 + 0.5*3,-0.4);

      \begin{scope}[shift={(7cm,0cm)}]
        \foreach \i in {0} {
          \draw[thick,fill=red] (\i*0.5,0) circle (0.2);
        };
        \foreach \i in {1,2,3,4,5,6,7,8,9,10} \draw[thick,fill=black!10] (\i*0.5,0) circle (0.2);
        \draw[thick] (-0.25,-0.3) rectangle (-0.25 + 0.5*5,0.3);
        \draw[thick] (-0.25 + 0.5*6,-0.3) rectangle (-0.25 + 0.5*11,0.3);
        \draw[->,thick] (0.2,0.4) to[bend left] (0.5*6-0.2,0.4);
      \end{scope}
    \end{tikzpicture}
  \end{center}

  The solution of this elementary combinatorial problem can be seen in the right picture. At least $\lfloor (n+1) /2 \rfloor$ adjacent transpositions are needed, and the minimum is obtained e.g. by choosing $i = 1$ and $j = \lfloor n/2 \rfloor$. The argument for $\SL(n,\bC)$ is the same except that the root spaces $\fg_\alpha$ are $2$--dimensional.
\end{Ex}

Similarly to the nonemptiness proof in \cite[Theorem 9.1]{GuichardWienhardDomains}, a bound on the dimension of the limit set can ensure that the domain is dense or connected.

\begin{Lem}\label{lem:dimensions}
  Let $\rho \colon \Gamma \to G$ be a $\Delta$--Anosov representation with limit map $\xi \colon \bdry \to \F_\Delta$. Then
  \begin{enumerate}
  \item $\dim\bdry \leq \mbcd(G)$.
  \item If $I \subset W_{\Delta,\eta}$ is a balanced ideal and $\dim\bdry \leq \mbcd(G) - 1$, then $\Omega(\xi(\bdry), I) \subset \F_\eta$ is dense.
  \item If $I \subset W_{\Delta,\eta}$ is a balanced ideal and $\dim\bdry \leq \mbcd(G) - 2$, then $\Omega(\xi(\bdry), I) \subset \F_\eta$ is connected.
  \end{enumerate}
\end{Lem}

\begin{Prf}
  We can assume that $\eta = \Delta$ in parts \itemnr{2} and \itemnr{3} as otherwise we could just lift to $\F_\Delta$. So let $I \subset W$ be a balanced ideal. We will calculate the covering dimension of
  \[\mathcal K = \F_\eta \setminus \Omega(\xi(\bdry), I) = \bigcup_{x \in \bdry} \bigcup_{w \in I} C_w(\xi(x)).\]
  Since $I$ is balanced, $\mathcal K$ is a continuous fiber bundle over $\bdry$ with fiber $\bigcup_{w \in I}C_w([1])$ (see \cite[Lemma B.8]{SteckerTreib} for details). Since the dimension can be calculated in local trivializations and the fiber is a CW--complex, $\dim \mathcal K = \dim\bdry + \dim \bigcup_{w\in I}C_w([1])$ \cite[Theorem 2]{Morita}. To bound the latter dimension, we use that $\dim C_w([1]) = \dim \Psi_w$ and that all $w \in I$ satisfy $w_0 w \not\leq w$. Furthermore, it is easy to check that $\dim \Psi_w = \dim \F_\Delta - \dim \Psi_{w_0w}$ for every $w \in W$. So we get the estimate
  \begin{equation}
    \dim \bigcup_{w\in I} C_w([1]) = \max_{w\in I} \dim \Psi_w \leq \max_{w_0 w \not\leq w} \dim \Psi_w = \dim\F_\Delta - \mbcd(G).
    \label{eq:dimension_estimate}
  \end{equation}

  Now if we assume $\dim\bdry \leq \mbcd(G) - 1$, then $\dim \mathcal K \leq \dim \F_\Delta - 1$. So $\Omega$ must be dense, as otherwise $\mathcal K$ would contain an open subset and therefore $\dim \mathcal K = \dim \F_\Delta$. This proves \itemnr{2}.

  For part \itemnr{3}, we can use Alexander duality \cite[Theorem 3.44]{Hatcher}: For a compact set $K$ of a closed manifold $M$, there is an isomorphism $H_i(M,M\setminus K;\bZ) \cong \check H^{n-i}(K;\bZ)$ for every $i$. Since $\dim \mathcal K \leq \dim \F_\Delta - 2$, and every Čech cohomology group above the covering dimension vanishes, we have $H_0(\F_\Delta, \Omega;\bZ) = \check H^n(\mathcal K;\bZ) = 0$ and $H_1(\F_\Delta, \Omega;\bZ) = \check H^{n-1}(\mathcal K;\bZ) = 0$. So by the long exact sequence of the pair $(\F_\Delta, \Omega)$ there is an isomorphsim $H_0(\Omega;\bZ) \cong H_0(\F_\Delta;\bZ)$, i.e. $\Omega$ is connected.

  Finally, for part \itemnr{1}, we just need a balanced ideal which gives equality in \eqref{eq:dimension_estimate}. This always exists: Take $w' \in W$ such that $w_0 w' \not\leq w'$ and which realizes the maximum. Then the ideal generated by $w'$ is slim and can therefore be extended to a balanced ideal $I$ by \cite[Lemma 3.34]{SteckerTreib} with $\max_{w \in I} \dim \Psi_w = \dim\Psi_{w'} = \max_{w_0w \not\leq w} \dim \Psi_w$. The corresponding $\mathcal K$ then satisfies $\dim\F_\Delta \geq \dim \mathcal K = \dim\bdry + \dim\F_\Delta - \mbcd(G)$, so $\dim\bdry \leq \mbcd(G)$.
\end{Prf}

\begin{Thm}\label{thm:cocompact_domains_correspondence}
  Let $\rho \colon \Gamma \to G$ be a $\Delta$--Anosov representation with limit map $\xi \colon \bdry \to \F_\Delta$ and $\Lambda = \xi(\bdry)$. Assume that $\dim \bdry \leq \mbcd(G) - 2$. Then every non--empty cocompact domain of discontinuity in $\F_\eta$ is dense and connected and there is a bijection
  \[\{\text{balanced ideals in $W_{\Delta,\eta}$}\} \leftrightarrow \{\text{non--empty cocompact domains of discontinuity in $\F_\eta$}\}\]
  given by $I \mapsto \Omega(\Lambda,I)$ and $\Omega \mapsto I(\Lambda, \Omega)$.
\end{Thm}

\begin{Prf}
  Let $\Omega$ be a non--empty cocompact domain of discontinuity. By \autoref{thm:cocompact_domains_lift} $\pi_\eta^{-1}(\Omega)$ is a union of connected components of $\Omega(\Lambda, \widetilde I)$ for some balanced ideal $\widetilde I \subset W$. But by \autoref{lem:dimensions} $\Omega(\Lambda,\widetilde I)$ is dense and connected. So $\pi_\eta^{-1}(\Omega) = \Omega(\Lambda, \widetilde I)$ and $\Omega$ is dense and connected.

  We have to prove that both maps are well--defined and inverses of each other. If $\Omega$ is a non--empty cocompact domain of discontinuity, then it is a union of connected components of some maximal domain $\widetilde \Omega \subset \F_\eta$ by \autoref{lem:cocompact_maximal}. Since $\Omega$ is dense it equals $\widetilde \Omega$ and is maximal itself. So by \autoref{pro:maximal_domains} $I(\Lambda, \Omega)$ is a fat ideal and $\Omega = \Omega(\Lambda, I(\Lambda,\Omega))$. Since $\Omega(\Lambda, I(\Lambda,\Omega))$ is dense and cocompact, \autoref{lem:not_cocompact} shows that $D(\Lambda, I(\Lambda,\Omega)) = \varnothing$, and by \autoref{lem:slim_no_intersections} this is equivalent to $I(\Lambda,\Omega)$ being slim, so $I(\Lambda,\Omega)$ is balanced.

  Conversely, if $I \subset W_{\Delta,\eta}$ is a balanced ideal, then $\Omega(\Lambda,I)$ is a cocompact domain of discontinuity by the main theorem of \cite{KapovichLeebPortiFlagManifolds}. It is dense and thus non--empty by \autoref{lem:dimensions}\itemnr{2}. By the above, $I(\Lambda,\Omega(\Lambda,I))$ is then a balanced ideal, and since $I \subset I(\Lambda,\Omega(\Lambda,I))$, they must be equal.
\end{Prf}

\subsection{A representation into \texorpdfstring{$\Sp(4,\bR)$}{Sp(4,R)} with infinitely many maximal domains}\label{sec:counterexample}

In this section, we describe an example of a representation (of a free group into $\Sp(4,\bR)$) which is Anosov (but not $\Delta$--Anosov) and where the analogue of \autoref{cor:maximal_domains2} does not hold, i.e. there are maximal domains of discontinuity which do not come from a balanced ideal. In fact, it will have infinitely many maximal domains of discontinuity, which are however not cocompact. It is unclear whether the cocompact domains for general Anosov representations can still be classified using balanced ideals.

Let $\Gamma = F_m$ be a free group in $m$ generators and $\rho_0 \colon \Gamma \to \SL(2,\bR)$ the holonomy of a compact hyperbolic surface with boundary. Such a representation is Anosov with a limit map $\xi_0 \colon \bdry \to \RP^1$ whose image $\Lambda_0 \coloneqq \xi_0(\bdry)$ is a Cantor set. We will now consider the representation $\rho = \iota \circ \rho_0$ into $\Sp(4,\bR)$, where
\[\iota \colon \SL(2,\bR) \to \Sp(4,\bR), \quad \begin{pmatrix}a & b \\ c & d\end{pmatrix} \mapsto \begin{pmatrix}a\mathbf{1} & b\mathbf{1} \\ c\mathbf{1} & d\mathbf{1}\end{pmatrix},\]
with $\mathbf{1}$ being the $2 \times 2$ identity matrix. Here we chose the symplectic form $\omega = \begin{psmallmatrix}0 & \mathbf{1} \\ -\mathbf{1} & 0\end{psmallmatrix}$. Then $\rho$ is $\{\alpha_2\}$--Anosov (where $\alpha_2$ is the simple root mapping a diagonal matrix to twice its lowest positive eigenvalue), but not $\Delta$--Anosov. Therefore, it carries a limit map $\xi \colon \bdry \to \Lag(\bR^4)$ to the manifold of Lagrangian subspaces.

The space $\Lag(\bR^4)$ admits the following (non--injective) parametrization:
\begin{align*}
  \Theta \colon \RP^1 \times \RP^1 \times \RP^1 &\to \Lag(\bR^4) \\
  \left(\left[\begin{smallmatrix}a \\ b\end{smallmatrix}\right], \left[\begin{smallmatrix}c \\ d\end{smallmatrix}\right], \left[\begin{smallmatrix}e \\ f\end{smallmatrix}\right]\right) & \mapsto \left\langle\begin{pmatrix}ae \\ af \\ be \\ bf\end{pmatrix}, \begin{pmatrix}cf \\ -ce \\ df \\ -de\end{pmatrix}\right\rangle.
\end{align*}
$\Theta$ is a smooth surjective map and the non--injectivity is precisely given by
\[\Theta(p,p,r) = \Theta(p,p,r') \quad \text{and} \quad \Theta(p,q,r) = \Theta(q,p,Rr) \qquad \forall p,q,r,r' \in \RP^1\]
with $R = \begin{psmallmatrix}0 & 1 \\ -1 & 0\end{psmallmatrix}$.


The action of $A \in \SL(2,\bR)$ on $\Lag(\bR^4)$ through $\iota$ is just
\[\iota(A) \Theta(p,q,r) = \Theta(Ap, Aq, r).\]
This results in the following simple description of the dynamical relations by the $\rho$--action on $\Lag(\bR^4)$:

\begin{Lem} \label{lem:Theta_dyn_rel}
  Let $\sim$ be the dynamical relation on $\Lag(\bR^4)$ by the action of $\rho(\Gamma)$. Then for all $x,y \in \Lambda_0$ and $p,q,r \in \RP^1$ we have
  \begin{equation}
    \label{eq:Theta_dyn_rel}
    \Theta(p,q,r) \sim \Theta(y,y,r) \quad \text{and} \quad \Theta(p, x, r) \sim \Theta(y, q, r).
  \end{equation}
  These are all dynamical relations.
\end{Lem}

\begin{Prf}
  Assume that $\Theta(p,q,r) \overset{(\rho(\gamma_n))}{\sim} \Theta(p',q',r')$ via a sequence $(\rho(\gamma_n)) \in \rho(\Gamma)^\bN$. Passing to subsequences, we can assume that $(\rho_0(\gamma_n)) \in \SL(2,\bR)^\bN$ is simply divergent with limits $(x,y) \in \Lambda_0$ and that there are sequences $(p_n), (q_n), (r_n) \in (\RP^1)^\bN$ such that
  \[p_n \to \widetilde p, \quad \rho_0(\gamma_n) p_n \to \widetilde p', \quad q_n \to \widetilde q, \quad \rho_0(\gamma_n) q_n \to \widetilde q', \quad r_n \to \widetilde r,\]
  \[\Theta(\widetilde p, \widetilde q, \widetilde r) = \Theta(p,q,r), \quad \Theta(\widetilde p', \widetilde q', \widetilde r) = \Theta(p',q',r').\]
  So $\widetilde p \sim \widetilde p'$ and $\widetilde q \sim \widetilde q'$ via $(\rho_0(\gamma_n))$. This either means that $\widetilde p = \widetilde q = x$ or $\widetilde p' = \widetilde q' = y$, in which case the relation is of the first type in \eqref{eq:Theta_dyn_rel}, or that $(\widetilde p, \widetilde q') = (x,y)$ or $(\widetilde q, \widetilde p') = (x,y)$, which is of the second type.

  Conversely, let $y \in \Lambda_0$ and $p,q \in \RP^1$. Since $|\Lambda_0| \geq 3$ we find $x \in \Lambda_0 \setminus \{p,q\}$, and by \autoref{lem:limit_pairs_in_Gamma} there is a sequence $(g_n) \in \rho_0(\Gamma)^\bN$ which is simply divergent with limits $(x,y)$. Then $p \overset{(g_n)}{\sim} y$ and $q \overset{(g_n)}{\sim} y$, which proves the first relation in \eqref{eq:Theta_dyn_rel}. For the second relation let $x,y \in \Lambda_0$ and $p,q,r \in \RP^1$. If $x = p$ or $y = q$ then it follows from the first relation. Otherwise, take a simply divergent sequence $(g_n) \in \rho_0(\Gamma)^\bN$ with limits $(x,y)$. Then $p \overset{(g_n)}{\sim} y$ and $q \overset{(g_n^{-1})}{\sim} x$. This shows the second relation in \eqref{eq:Theta_dyn_rel}.
\end{Prf}

\begin{Prop}
  Let $A \subset \RP^1$ be a minimal closed subset such that $A \cup RA = \RP^1$. Then
  \[\Omega_A = \Lag(\bR^4) \setminus \{\Theta(p,q,r) \mid p \in \Lambda_0, q \in \RP^1, r \in A\}\]
  is a maximal domain of discontinuity for $\rho$.
\end{Prop}

\begin{Prf}
  $\Omega_A$ is open since $(\RP^1)^3$ is compact and $\Theta$ therefore is a closed map.

  Assume that there was a dynamical relation within $\Omega_A$. Then following \eqref{eq:Theta_dyn_rel} it would either be of the form $\Theta(p,q,r) \sim \Theta(x,x,r)$ or $\Theta(p,x,r) \sim \Theta(y,q,r)$ with $x,y \in \Lambda_0$ and $p,q,r \in \RP^1$. In the first case, $\Theta(x,x,r)$ is independent of $r$, so we can assume $r \in A$, and $\Theta(x,x,r)$ can thus not be in $\Omega_A$. In the second case, $\Theta(p,x,r) = \Theta(x,p,Rr)$ can be in $\Omega_A$ only if $r \not\in RA$ and $\Theta(y,q,r) \in \Omega_A$ implies $r \not\in A$. But by assumption both can not hold at the same time.

  Finally, assume that $\Omega_A$ was not maximal, i.e. there was another domain of discontinuity $\Omega \subset \Lag(\bR^4)$ with $\Omega_A \subsetneq \Omega$. Let
  \[A' = \{r \in \RP^1 \mid \forall x \in \Lambda_0, q \in \RP^1 \colon \Theta(x,q,r) \not \in \Omega\}.\]
  Then $A' \subsetneq A$ and $A'$ is closed. Since $A$ is minimal among closed sets with $A \cup RA = \RP^1$, there has to exist some $r \in \RP^1 \setminus (A' \cup RA')$. But since $r, Rr \not\in A'$ then there are $x,y \in \Lambda_0$ and $p,q \in \RP^1$ with $\Theta(y,q,r), \Theta(x,p,Rr) \in \Omega$. But these are dynamically related by \autoref{lem:Theta_dyn_rel}, a contradiction.
\end{Prf}

Through the accidental isomorphism $\PSp(4,\bR) \cong \SO_0(2,3)$ the space $\Lag(\bR^4)$ can be identified with the space of isotropic lines in $\bR^{2,3}$. The form of signature $(2,3)$ restricts to a Lorentzian metric on this space, which is why it is also called the $(2+1)$ Einstein universe. A detailed explanation of its geometry can be found in \cite[Section 5]{BarbotCharetteDrummGoldmanMelnick}.

We can use this to visualize $\Theta$ and the construction of $\Omega_A$ above: The limit set $\xi(\bdry) \subset \Lag(\bR^4)$ is a Cantor set on the line $\{\Theta(x,x,*) \mid x \in \RP^1\}$. If we take two different points on this line, described by $x,y \in \RP^1$, their light cones intersect in the circle $\{\Theta(x,y,r) \mid r \in \RP^1\}$, where $r$ acts as a global angle coordinate. If $x,y \in \Lambda_0$ then every point on the future pointing light ray emanating from $x$ in a direction $r$ is dynamically related to every point on the past pointing light ray from $y$ in direction $r$ (the red and blue lines in \autoref{fig:einstein_universe}). So by choosing the set $A \subset \RP^1$, we decide for every angle whether to take out from our domain all the future or all the past pointing light rays emanating from the points in the limit set in this direction.

\begin{center}
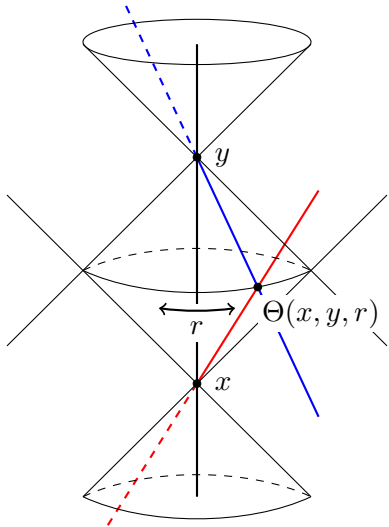

  \begin{tikzpicture}
    \draw[thick,red] (0,-1.5) -- (0.8*2,-1.5 + 1.28*2);
    \draw[thick,red,dashed] (0,-1.5) -- (-0.8*1.5,-1.5 - 1.28*1.5);
    \draw[thick,blue] (0,1.5) -- (0.8*2,1.5 - 1.72*2);
    \draw[thick,blue,dashed] (0,1.5) -- ({0.8*(-1.2)},{1.5 - 1.72*(-1.2)});

    \draw[thick] (0,3) -- (0,-3);
    \draw[fill=black] (0,1.5) circle (0.05) node[right,xshift=1mm] {$y$};
    \draw[fill=black] (0,-1.5) circle (0.05) node[right,xshift=1mm] {$x$};
    \draw (-1.5,-3) -- (2.5, 1);
    \draw (1.5,-3) -- (-2.5, 1);
    \draw (-1.5,3) -- (2.5, -1);
    \draw (1.5,3) -- (-2.5, -1);

    \draw[fill=black] (0.8,-0.22) circle (0.05) node[anchor=north west,fill=white,inner sep=0.5mm,yshift=-1mm] {$\Theta(x,y,r)$};

    \draw (-1.5,0) arc [start angle = -135, end angle = -45, x radius = 3/(cos(45)-cos(135)), y radius = 1];
    \draw[dashed] (-1.5,0) arc [start angle = 135, end angle = 45, x radius = 3/(cos(45)-cos(135)), y radius = 1];

    \draw (0,3.03) ellipse [x radius = 1.5, y radius = 0.3];

    \draw (-1.5,-3) arc [start angle = -135, end angle = -45, x radius = 3/(cos(45)-cos(135)), y radius = 1];
    \draw[dashed] (-1.5,-3) arc [start angle = 135, end angle = 45, x radius = 3/(cos(45)-cos(135)), y radius = 1];

    \fill[white] (-0.2,-0.45) rectangle (0.2,-0.95);
    \draw[<->,thick] (-0.5,-0.5) arc [start angle = -105, end angle = -75, x radius = 1/(cos(75)-cos(105)), y radius = 1] node[midway,below] {$r$};
  \end{tikzpicture}
  \captionof{figure}{The parametrization $\Theta$ interpreted by intersecting light cones in $\Lag(\bR^4)$. The vertical line is the set of points $\Theta(x,x,*)$ containing the limit set.}
  \label{fig:einstein_universe}
\end{center}

\section{\texorpdfstring{Representations into $\SL(n,\bR)$ or $\SL(n,\bC)$}{Representations into SL(n,R) or SL(n,C)}}

\subsection{Balanced ideals}\label{sec:ideals_in_SL}

The question for which $\eta \subset \Delta$ there exists a balanced ideal in $W_{\Delta,\eta}$ is only combinatorial. For $G = \SL(n,\bK)$ with $\bK \in \{\bR, \bC\}$ the answer is given by the following proposition. \Autoref{thm:intro_Hitchin} and its corollaries then immediately follow using \autoref{thm:cocompact_domains_correspondence}.

A maximal compact subgroup of $\SL(n,\bR)$ is $K = \SO(n)$ and for $G = \SL(n,\bC)$ we can choose $K = \SU(n)$. In either case, a maximal abelian subalgebra $\fa$ of $\fso(n)^\perp$ resp. $\fsu(n)^\perp$ are the traceless real diagonal matrices, and a simple system of restricted roots is given by $\{\alpha_i = \lambda_i - \lambda_{i+1}\}$, where $\lambda_i \colon \fa \to \bR$ maps to the $i$--th diagonal entry.

\begin{Prop}\label{pro:SLn_domains}
  Let $\eta = \{\alpha_{i_1}, \dots, \alpha_{i_k}\} \subset \Delta$ be a subset of the simple roots of $\SL(n,\bK)$, with $0 = i_0 < i_1 < \dots < i_k < i_{k+1} = n$. Let
  \[\delta = |\{0 \leq j \leq k \mid i_{j+1} - i_j \;\text{\upshape is odd}\}|.\]
  If $n$ is even, a balanced ideal exists in $W_{\Delta,\eta}$ if and only if $\delta \geq 1$. If $n$ is odd, a balanced ideal exists in $W_{\Delta,\eta}$ if and only if $\delta \geq 2$.
\end{Prop}

\begin{Prf}
  A balanced ideal exists if and only if the action of $w_0$ on $W_{\Delta,\eta}$ by left--mul\-ti\-pli\-ca\-tion fixes no element of $W_{\Delta,\eta}$ (see \cite[Proposition 3.29]{KapovichLeebPortiFlagManifolds} or \cite[Lemma 3.34]{SteckerTreib}). This means that $ww_0w^{-1} \not\in \langle \ceta \rangle$ for any $w \in W$. The Weyl group of $\SL(n,\bK)$ can be identified with the symmetric group $S_n$ with its generators $\Delta$ being the adjacent transpositions. Assume first that $n$ is even. Then $w_0$ is the order--reversing permutation and its conjugates are precisely the fixed point free involutions in $S_n$. So the existence of balanced ideals is equivalent to every involution in $\langle \ceta \rangle$ having a fixed point.

  Now observe that $\langle \ceta \rangle$ is a product of symmetric groups, namely $\langle \ceta \rangle \! \cong \prod_{j = 0}^k S_{i_{j+1} - i_j}$, and that there are fixed point free involutions in $S_k$ if and only if $k$ is even. So a balanced ideal exists iff at least one of the $i_{j+1} - i_j$ is odd, i.e. $\delta \geq 1$.

  The same argument works if $n$ is odd, except that the conjugates of $w_0$ are then involutions in $S_n$ with precisely one fixed point (every involution has at least one), and so we need $\delta \geq 2$ to have none of these in $\langle \ceta \rangle$.
\end{Prf}

For the action on Grassmannians, \autoref{pro:SLn_domains} specializes to the following simple condition: A balanced ideal exists in $W_{\Delta,\{\alpha_k\}}$ if and only if $n$ is even and $k$ is odd. We can enumerate all balanced ideals for $n \leq 10$ using a computer and obtain the following number of balanced ideals in $W_{\Delta,\{\alpha_k\}}$:

\begin{center}
  \begin{tabular}{l|ccccc}
    & $k=1$ & $k=3$ & $k=5$ & $k=7$ & $k=9$ \\
    \hline
    $n=2$ & $1$  &&&& \\
    $n=4$ & $1$ & $1$ &&& \\
    $n=6$ & $1$ & $2$ & $1$ && \\
    $n=8$ & $1$ & $7$ & $7$ & $1$ &\\
    $n=10$ & $1$ & $42$ & $2227$ & $42$ & $1$ \\
  \end{tabular}
  \label{fig:number_balanced_ideals}
\end{center}

In particular, a cocompact domain of discontinuity in projective space $\RP^{n-1}$ or $\CP^{n-1}$ exists if and only if $n$ is even. Interestingly, these are also precisely the dimensions which admit complex Schottky groups by \cite{Cano}.

\subsection{Hitchin representations}\label{sec:Hitchin}

Let $\Gamma$ be the fundamental group of a closed surface. A \emph{Hitchin representation} $\rho \colon \Gamma \to \SL(n,\bR)$ is a representation which can be continuously deformed to a representation of the form $\iota \circ \rho_0$ where $\rho_0 \colon \Gamma \to \SL(2,\bR)$ is discrete and injective and $\iota \colon \SL(2,\bR) \to \SL(n,\bR)$ is the irreducible representation. Hitchin representations are $\Delta$--Anosov \cite{Labourie}.

\Autoref{thm:cocompact_domains_correspondence} together with \autoref{ex:mbcd_SLn} shows that if $n \geq 5$ then the cocompact domains of discontinuity of a Hitchin representation in any flag manifold $\F_\eta$ are in 1:1 correspondence with the balanced ideals in $W_{\Delta,\eta}$. These were discussed in \autoref{sec:ideals_in_SL}.

For completeness, let us also have a look at the cases $n \in \{2,3,4\}$.

In $\SL(2,\bR)$ the Hitchin representations are just the discrete injective representations. The only flag manifold is $\RP^1$, and since the limit maps $\xi \colon \bdry \to \RP^1$ of Hitchin representations are homeomorphisms, there can be no non--empty domain of discontinuity in $\RP^1$.

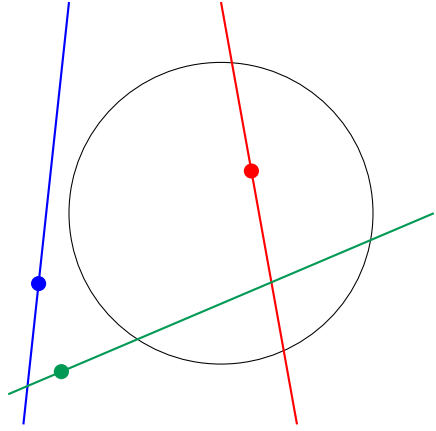
\begin{figure}
  \centering
  \begin{tikzpicture}[scale=2]
    \draw (0,0) circle (1.0);
    \draw[thick,red] (0,1.4) -- (0.5,-1.4);
    \fill[red] (0.2,0.28) circle (0.05);
    \draw[thick,blue] (-1,1.4) -- (-1.3,-1.4);
    \fill[blue] (-1.2,-0.4666) circle (0.05);
    \draw[thick,ForestGreen] (1.4,0) -- (-1.4, -1.2);
    \fill[ForestGreen] (-1.05,-1.05) circle (0.05);
  \end{tikzpicture}
  \caption{The three connected components of the maximal domain of discontinuity in $\F_\Delta$ for a Hitchin representation $\rho \colon \Gamma \to \SL(3,\bR)$. One exemplary flag out of every component is shown, appearing as a point on a line in $\RP^2$.}
  \label{fig:SL3_components}
\end{figure}

In the case of $\SL(3,\bR)$ there is only a single balanced ideal $I \subset W$. By \autoref{thm:cocompact_domains_lift} the lift of any cocompact domain of discontinuity to the full flag manifold $\F_\Delta$ must be a union of connected components of the corresponding domain, which is
\[\Omega(\Lambda,I) = \{f \in \F_\Delta \mid \forall x \in \bdry \colon f^{1} \neq \xi^1(x) \land f^{2} \neq \xi^2(x)\}.\]
It is known (see \cite{ChoiGoldmanClosed}) that for any Hitchin representation $\rho$ into $\SL(3,\bR)$ there exists a properly convex open domain $D \subset \RP^2$ on which $\rho$ acts properly discontinuously and cocompactly. The image of the limit map of $\rho$ are then the flags consisting of a point on $\partial D$ and the tangent line of $D$ through this point. The domain $\Omega(\Lambda, I) \subset \F_\Delta$ therefore splits into three connected components: One of them (shown in red in \autoref{fig:SL3_components}) consists of flags (i.e. a point and a line through it in $\RP^2$) with the point inside $D$. The second component (blue in \autoref{fig:SL3_components}) are flags whose line avoids $D$, and the third (green) consists of flags whose line goes through $D$ but with the point being outside.

Only the red component descends to a domain in $\RP^2$, and only the blue one to $\Gr(2,3)$, each forming the unique cocompact domain of discontinuity in these manifolds. The cocompact domains in the full flag manifold are any unions of one or more of the three components.

Finally, let's have a look at $\SL(4,\bR)$. There are ten balanced ideals in $W$ in this case (see \autoref{sec:appendix_A3}), corresponding to ten maximal domains of discontinuity in $\F_\Delta$. By dimension arguments as in the proof of \autoref{lem:dimensions} all of these domains are dense and eight of them are connected. The other two domains are
\[\Omega_1 = \{f \in \F_\Delta \mid \forall x \in \bdry \colon f^1 \not\subset \xi^2(x)\}, \quad \Omega_2 = \{f \in \F_\Delta \mid \forall x \in \bdry \colon \xi^2(x) \not\subset f^3\}.\]
The topology of these domains does not change when the representation is continuously deformed, and from the Fuchsian case we can easily see that $\Omega_1$ and $\Omega_2$ each have two connected components, all of which are lifts of domains in $\RP^3$ or $\Gr(3,4)$, respectively. The quotient of one of the components in $\RP^3$ describes a convex foliated projective structure on the unit tangent bundle of $S$ \cite{GuichardWienhardConvexFoliated}.

Out of the other 8 cocompact domains in $\F_\Delta$, one descends to the partial flag manifolds $\F_{1,2}$, $\F_{2,3}$ and $\F_{1,3}$, each. Counting all possible combinations of connected components separately, we have $14$ different non--empty cocompact domains of discontinuity in $\F_\Delta$.

\section{\texorpdfstring{Symplectic Anosov representations into $\Sp(4n+2,\bR)$}{Symplectic Anosov representations into Sp(4n+2,R)}}

In this section, we will prove \autoref{thm:intro_compactification}, i.e. construct a compactification for locally symmetric spaces modeled on the bounded symmetric domain compactification. We first recall some facts on this compactification, which can be found in \cite{Helgason} and \cite{Satake}.


Every Hermitian symmetric space can be realized as a bounded symmetric domain in some $\bC^N$. That is an open, connected and bounded subset $D \subset \bC^N$ such that for every point $x \in D$ there is an involutive holomorphic diffeomorphism from $D$ to itself which has $x$ as an isolated fixed point. Concretely, to get the symmetric space $\Sp(2n,\bR)/\U(n)$ we can consider the bounded symmetric domain
\[D = \{Z \in \Sym(n, \bC) \mid 1 - \overline Z Z \;\text{is positive definite}\} \subset \bC^{n(n+1)/2}.\]
The group of holomorphic diffeomorphisms of $D$ is isomorphic to $\Sp(2n,\bR)$ and acts with stabilizer $\U(n)$. We compactify the symmetric space by taking the closure $\overline D$ in $\bC^{n(n+1)/2}$. This is the \emph{bounded symmetric domain compactification} of $\Sp(2n,\bR)/\U(n)$.

Instead of working with bounded symmetric domains, we will use an equivalent model of $\Sp(2n,\bR)/\U(n)$, the \emph{Siegel space}. Let $\omega$ be a symplectic form on $\bR^{2n}$ and $\omega_\bC$ its complexification on $\bC^{2n}$. Together with the real structure this defines an (indefinite) Hermitian form
\[h(v,w) \coloneqq \frac{i}{2}\omega_\bC(\overline v, w) \qquad \forall v, w \in \bC^{2n}.\]
The Siegel space is the subspace $\mathcal H_{n0} \subset \Lag(\bC^{2n})$ of complex Lagrangians $L$ such that $h|_{L\times L}$ is positive definite.

The correspondence between these models uses the Cayley transform: Regard the symmetric complex matrices $\Sym(n,\bC)$ as a subset of $\Lag(\bC^{2n})$ by mapping $X \in \Sym(n,\bC)$ to $\{(Xv,v) \mid v \in \bC^n\} \in \Lag(\bC^{2n})$. The Cayley transform on $\Lag(\bC^{2n})$ is just the action of the matrix
\[\frac{e^{i\pi/4}}{\sqrt{2}}\begin{pmatrix} -i & i \\ 1 & 1\end{pmatrix} \in \Sp(2n,\bC).\]
It maps $D$ to $\mathcal H_{n0}$ and $\overline D$ to $\overline{\mathcal H_{n0}}$, establishing an equivalence of these compactifications.

More generally, let $\mathcal H_{pq} \subset \Lag(\bC^{2n})$ be the set of Lagrangians such that $h$ restricted to them has signature $(p,q)$, meaning that there is an orthogonal basis with $p$ vectors of positive norm and $q$ vectors of negative norm (and possibly null vectors). Then
\[\Lag(\bC^{2n}) = \bigsqcup_{\stackrel{0 \leq p,q \leq n}{p+q \leq n}} \mathcal H_{pq},\]
and the $\mathcal H_{pq}$ are precisely the orbits of the action of $\Sp(2n,\bR) \subset \Sp(2n,\bC)$ on $\Lag(\bC^{2n})$. Furthermore, the map
\[\mathcal H_{pq} \to \Is_{n-p-q}(\bR^{2n}), \quad L \mapsto L \cap \overline L\]
makes every $\mathcal H_{pq}$ a fiber bundle over the isotropic $(n\mathord{-}p\mathord{-}q)$--subspaces with fiber the semi--Riemannian symmetric space $\Sp(2p + 2q,\bR) / \U(p,q)$. In particular, this means $\mathcal H_{n0} = \Sp(2n,\bR)/\U(n)$. A more detailed explanation can be found in \cite{WienhardRepresentationsGeometricStructures}.

Now let $\rho \colon \Gamma \to \Sp(2n,\bR)$ be an $\{\alpha_{n}\}$--Anosov representation and $\xi \colon \Gamma \to \F_{\{\alpha_n\}} = \Lag(\bR^{2n})$ its limit map. Important examples of these representations are maximal representations from a surface group $\Gamma = \pi_1S$, if either $S$ is a closed surface, or an open surface and the boundary elements map to Shilov hyperbolic elements of $\Sp(2n,\bR)$. The following theorem implies \autoref{thm:intro_compactification}:


\begin{Thm}\label{thm:compactification}
  If $n$ is odd, then there exists a balanced ideal $I \subset W_{\{\alpha_n\},\{\alpha_n\}}$. Therefore,
  \[\widehat X \coloneqq \overline {\mathcal H_{n0}} \cap \Omega(\xi(\bdry),I) = \{ L \in \overline{\mathcal H_{n0}} \mid \dim_\bC L \cap \xi(x)^\bC < n/2 \; \forall x \in \bdry\}\]
  is $\Gamma$--invariant and the quotient $\Gamma \backslash \widehat X$ is a compactification of the locally symmetric space $\Gamma \backslash \mathcal H_{n0} = \Gamma \backslash \Sp(2n,\bR) / \U(n)$.
\end{Thm}

\begin{Prf}
  The Weyl group $W$ of $\Sp(2n,\bC)$ can be identified with the group of permutations $\pi$ of $\{-n, \dots, -1,1, \dots, n\}$ with $\pi(-i) = -\pi(i)$ for all $i$. In this identification, the generator $\alpha_n$ negates $n$ and $-n$ keeping everything else fixed, while $\alpha_k$ for $k \neq n$ exchanges $k$ with $k+1$ and $-k$ with $-k-1$. The longest element $w_0$ negates everything.

  Denote by $[k] \in W_{\{\alpha_n\},\{\alpha_n\}} = \langle \alpha_1,\dots,\alpha_{n-1} \rangle \backslash W / \langle \alpha_1,\dots,\alpha_{n-1} \rangle$ the equivalence class of permutations which map exactly $k$ positive numbers to positive ones. Then $w_0 [k] = [n-k]$ and $[k] \leq [\ell]$ in the Bruhat order if and only if $k \geq \ell$. Furthermore, $\pos(L,L') = [k]$ for two Lagrangians $L,L' \in \Lag(\bC^{2n})$ if and only if $\dim_\bC (L \cap L') = k$.

  If $n$ is odd, then $I = \{[k] \mid k > n/2\}$ is a balanced ideal, so
  \[\Omega \coloneqq \Omega(\xi(\bdry),I) = \{L \in \Lag(\bC^{2n}) \mid \dim_\bC(L \cap \xi(x)^\bC) < n/2 \;\forall x \in \bdry\}\]
  is a cocompact domain of discontinuity for $\rho$ by \cite{KapovichLeebPortiFlagManifolds}. Here $\rho$ is regarded as a representation into $\Sp(2n,\bC)$. But because $\rho$ maps into $\Sp(2n,\bR)$, it preserves $\mathcal H_{pq}$ and therefore also $\widehat X$. The quotient $\Gamma \backslash \widehat X$ is a closed subset of $\Gamma \backslash \Omega$ and thus also compact.
\end{Prf}

Note that this is just a special case of a general principle to construct compactifications for locally symmetric spaces arising from Anosov representations described in \cite{GuichardKasselWienhard} and \cite{KapovichLeebFinsler}. However, it is a particularly interesting one as this compactification is modeled on the bounded symmetric domain compactification for $\Sp(2n,\bR)/\U(n)$.

\appendix

\section{Lists of balanced ideals}\label{sec:list_of_balanced_ideals}

The description of cocompact domains using balanced ideals often reduces finding them to enumerating balanced ideals, which is a purely combinatorial problem. This means we can use a computer to do it. This section shows the resulting lists and numbers in some potentially interesting cases. The program used to compute them was written by the author together with David Dumas. It can be found online at \url{https://florianstecker.de/balancedideals/}.

\subsection{Balanced ideals in \texorpdfstring{$A_n$}{An}}

Assume the Coxeter system $(W, \Delta)$ defined by the Weyl group and the restricted roots of $G$ is of type $A_n$. For example, $G$ could be the group $\SL(n+1,\bR)$ or $\SL(n+1,\bC)$. We write $\Delta = \{\alpha_1,\dots,\alpha_n\}$ in such a way that $\alpha_i\alpha_{i+1} \in W$ is an element of order $3$. The following tables show all balanced ideals in $W$. For every balanced ideal $I$, it shows the subset of $\Delta$ it is left-- and right--invariant by. This means that for $\theta, \eta \subset \Delta$ the balanced ideals in $W_{\theta,\eta}$ are precisely those in $W$ whose left--invariance includes $\ctheta$ and right--invariance includes $\ceta$. We also show the (real resp. complex) dimension of the set we have to take out of the $\F_\Delta$ for every limit point, and a minimal set of elements of $W$ generating $I$ as an ideal.

\subsubsection{Balanced ideals in \texorpdfstring{$A_1$}{A1}}

\begin{center}
  \begin{tabular}{c|c|c|l}
    \textbf{left--invariance} & \textbf{right--invariance} & \textbf{dimension} & \textbf{generators} \\
    \hline
    $\varnothing$ & $\varnothing$ & 0 & $1$ \\
  \end{tabular}
\end{center}

\subsubsection{Balanced ideals in \texorpdfstring{$A_2$}{A2}}

\begin{center}
  \begin{tabular}{c|c|c|l}
    \textbf{left--invariance} & \textbf{right--invariance} & \textbf{dimension} & \textbf{generators} \\
    \hline
    $\varnothing$ & $\varnothing$ & 1 & $\alpha_1, \alpha_2$ \\
  \end{tabular}
\end{center}

\subsubsection{Balanced ideals in \texorpdfstring{$A_3$}{A3}}\label{sec:appendix_A3}

\begin{center}
  \begin{tabular}{c|c|c|l}
    \textbf{left--invariance} & \textbf{right--invariance} & \textbf{dimension} & \textbf{generators} \\
    \hline
    $\{\alpha_1,\alpha_3\}$ & $\{\alpha_1,\alpha_2\}$ & 4 & $\alpha_3 \alpha_1 \alpha_2 \alpha_1$ \\
    $\{\alpha_1,\alpha_3\}$ & $\{\alpha_2,\alpha_3\}$ & 4 & $\alpha_1 \alpha_2 \alpha_3 \alpha_2$ \\
    $\{\alpha_2\}$ & $\varnothing$ & 3 & $\alpha_1\alpha_2\alpha_1,\; \alpha_2\alpha_1\alpha_3,\; \alpha_2\alpha_3\alpha_2$ \\
    $\varnothing$ & $\{\alpha_1\}$ & 3 & $\alpha_3\alpha_2\alpha_1,\; \alpha_1\alpha_2\alpha_1,\; \alpha_2\alpha_1\alpha_3$ \\
    $\varnothing$ & $\{\alpha_2\}$ & 3 & $\alpha_1\alpha_3\alpha_2,\; \alpha_1\alpha_2\alpha_1,\; \alpha_2\alpha_3\alpha_2$ \\
    $\varnothing$ & $\{\alpha_3\}$ & 3 & $\alpha_1\alpha_2\alpha_3,\; \alpha_2\alpha_1\alpha_3,\; \alpha_2\alpha_3\alpha_2$ \\
    $\varnothing$ & $\varnothing$ & 3 & $\alpha_3\alpha_2\alpha_1,\; \alpha_1\alpha_3\alpha_2,\; \alpha_1\alpha_2\alpha_3$ \\
    $\varnothing$ & $\varnothing$ & 3 & $\alpha_3\alpha_2\alpha_1,\; \alpha_1\alpha_2\alpha_3,\; \alpha_2\alpha_1\alpha_3$ \\
    $\varnothing$ & $\varnothing$ & 3 & $\alpha_3\alpha_2\alpha_1,\; \alpha_1\alpha_3\alpha_2,\; \alpha_1\alpha_2\alpha_1,\; \alpha_2\alpha_3$ \\
    $\varnothing$ & $\varnothing$ & 3 & $\alpha_1\alpha_3\alpha_2,\; \alpha_1\alpha_2\alpha_3,\; \alpha_2\alpha_3\alpha_2,\; \alpha_2\alpha_1$
  \end{tabular}
\end{center}

\subsubsection{The number of balanced ideals in \texorpdfstring{$A_4$}{A4}}

There are $4608$ balanced ideals in $W$, so we cannot list them all. Instead, the following table shows just how many balanced ideals exist in $W_{\theta,\eta}$ for any choice of $\theta, \eta \subset \Delta$ with $\iota(\theta) = \theta$. The rows correspond to different values of $\theta$ (for example \textbf{14} stands for $\theta = \{\alpha_1,\alpha_4\}$) while the columns correspond to $\eta$.

\begin{center}
  \begin{tabular}{l|rrrrrrrrrrrrrrr}
    & \textbf{1234} & \textbf{123} & \textbf{134} & \textbf{124} & \textbf{234} & \textbf{12} & \textbf{13} & \textbf{14} & \textbf{23} & \textbf{24} & \textbf{34} & \textbf{1} & \textbf{2} & \textbf{3} & \textbf{4} \\
    \hline
    \textbf{1234} & 4608 & 35 & 57 & 57 & 35 & 2 & 0 & 3 & 0 & 0 & 2 & 0 & 0 & 0 & 0 \\
    \textbf{14} & 1  & 0 & 0 & 0 & 0 & 0 & 0 & 0 & 0 & 0 & 0 & 0 & 0 & 0 & 0 \\
    \textbf{23} & 12 & 2 & 5 & 5 & 2 & 1 & 0 & 2 & 0 & 0 & 1 & 0 & 0 & 0 & 0
  \end{tabular}
\end{center}

One feature stands out: There is only a single balanced ideal in $W_{\{\alpha_1,\alpha_4\},\Delta}$, and it has no right--invariances at all. In fact, we have the same situation generally in $W_{\{\alpha_1,\alpha_n\},\Delta}$ if $W$ is of type $A_n$. For an $\{\alpha_1,\alpha_n\}$--Anosov representation into $\SL(n+1,\bR)$ with limit map $\xi \colon \bdry \to \F_{1,n}$, this corresponds to the cocompact domain
\[\Omega = \F_\Delta \setminus \{F \in \F_\Delta \mid \exists x \in \bdry, i \leq n \colon \xi^{(1)}(x) \subset F^{(i)} \subset \xi^{(n)}(x)\}\]
which was also constructed in \cite[10.2.3]{GuichardWienhardDomains} using the adjoint representation.

\subsection{\texorpdfstring{$\{\alpha_1,\dots,\alpha_{p-1}\}$--Anosov representations into $\SO_0(p,q)$}{\{alpha\_1,...,alpha\_p-1\}-Anosov representations into SO(p,q)}}

Guichard and Wienhard recently identified an interesting class of surface group representations they call $\Theta$--positive representations \cite{GuichardWienhardPositivity}. This includes a family of representations into $\SO_0(p,q)$ with $p < q$ which they conjecture to be a union of connected components and to be $\theta$--Anosov with $\theta = \{\alpha_1,\dots,\alpha_{p-1}\}$. Here we ordered the simple roots such that non--consecutive ones commute and $\alpha_{p-1}\alpha_p$ has order 4. If this conjecture is true, balanced ideals in $W_{\theta,\eta}$ induce cocompact domains of discontinuity of these representations. Similarly to the table in \autoref{sec:ideals_in_SL}, the following table shows the number of balanced ideals in $W_{\theta,\eta}$ for $\eta = \{\alpha_k\}$, i.e. corresponding to domains in Grassmannians $\Is_k(\bR^{p,q})$ of isotropic $k$--subspaces.

\begin{center}
  \begin{tabular}{r|rrrrrrr}
    & $k = 1$ & $k = 2$ & $k = 3$ & $k = 4$ & $k = 5$ & $k = 6$ & $k = 7$ \\
    \hline
    $p = 2$ & 0 & 1 &&&&& \\
    $p = 3$ & 0 & 1 & 1 &&&& \\
    $p = 4$ & 0 & 1 & 2 & 2 &&& \\
    $p = 5$ & 0 & 1 & 7 & 14 & 3 && \\
    $p = 6$ & 0 & 1 & 42 & 616 & 131 & 7 & \\
    $p = 7$ & 0 & 1 & 429 & 303742 & 853168 & 8137 & 21 \\
  \end{tabular}
\end{center}

There is always a unique balanced ideal for $\eta = \{\alpha_2\}$. It corresponds to the cocompact domain of discontinuity
\[\Omega = \{V \in \Is_2(\bR^{p,q}) \mid V \perp \xi^{(i)}(x) \Rightarrow V \cap \xi^{(i)}(x) = 0 \ \forall x \in \bdry\; \forall i \leq p-1\}\]
in the space of isotropic planes.

\subsection{\texorpdfstring{$\{\alpha_1,\alpha_2\}$--Anosov representations into $F_4$}{F\_4}}

There is another exceptional family of $\Theta$--positive representations, which are conjectured to be $\{\alpha_1,\alpha_2\}$--Anosov in a group $G$ with Weyl group of type $F_4$. This table shows the number of balanced ideals in $W_{\{\alpha_1,\alpha_2\},\eta}$ for different choices of $\eta$. Again, \textbf{134} is a shorthand for $\eta = \{\alpha_1,\alpha_3,\alpha_4\}$.

\begin{center}
  \begin{tabular}{l|rrrrrrrrrrrrrrr}
    & \textbf{1234} & \textbf{123} & \textbf{134} & \textbf{124} & \textbf{234} & \textbf{12} & \textbf{13} & \textbf{14} & \textbf{23} & \textbf{24} & \textbf{34} & \textbf{1} & \textbf{2} & \textbf{3} & \textbf{4} \\
    \hline
    \textbf{12} & 1270 & 182 & 140 & 66 & 44 & 16 & 18 & 5 & 14 & 6 & 4 & 1 & 2 & 2 & 0 \\
  \end{tabular}
\end{center}

\KOMAoption{fontsize}{10pt}
\printbibliography

\vspace{1cm}

Florian Stecker \\
Heidelberg Institute for Theoretical Studies \\
Mathematics Department, Heidelberg University \\
\emph{E--mail address:} \texttt{fstecker@mathi.uni-heidelberg.de}

\end{document}